\documentclass[11pt]{article}
\usepackage{amsfonts}
\usepackage{epsfig}
\usepackage{amsmath}

\newcommand{\beq}[1]{\begin{equation}\label{#1}}
\newcommand{\eeq}{\end{equation}}
\newcommand {\calleq}[1]{(\ref{#1})}
\newcommand{\Caption}[2]{\caption{\it #2\hfill \label{#1}}}
\newcommand{\bl}{\vspace{1mm}}
\let\eps=\varepsilon
\let\phi=\varphi

\let\Omg=\Omega
\newcommand\hchi{\lambda}
\newcommand\chib{\overline{\lambda}}
\newcommand\ctoda{\mu}
\newcommand\ts{h}
\newcommand\cO{{\cal O}}
\newcommand{\cV}{{\cal V}}
\newcommand{\ab}{$\alpha+\beta$}
\newcommand{\bt}{$\beta_T$}
\newcommand{\gt}{$\gamma_T$}
\newcommand{\av}{variable $\alpha$}
\newcommand{\pb}{pure $\beta$}
\newcommand{\gd}{$\gamma+\delta$}
\newcommand{\pd}{pure $\delta$}
\newcommand\dist{\mathrm{dist}\,}
\newcommand\media[1]{\langle #1 \rangle}
\begin{document}

\title{\bf The Fermi--Pasta--Ulam problem \\ 
and its underlying integrable dynamics:\\
an approach through Lyapunov Exponents}

\author{
G.~Benettin\thanks{Universit\`a di Padova, Dipartimento di Matematica
  ``Tullio Levi-Civita'', Via Trieste 63, 35121 Padova (Italy);
    {\tt benettin@math.unipd.it}}, 
S.~Pasquali\thanks{Universit\`a di Roma Tre, Dipartimento di
  Matematica e Fisica; {\tt spasquali@mat.uniroma3.it}},
A.~Ponno\thanks{Universit\`a di Padova, Dipartimento di Matematica
  ``Tullio Levi-Civita'', Via Trieste 63, 35121 Padova (Italy);
  {\tt ponno@math.unipd.it}}
} 

\date{}

\maketitle
\vskip 3truecm
\begin{abstract}
FPU models, in dimension one, are perturbations either of the linear
model or of the Toda model; perturbations of the linear model include
the usual $\beta$-model, perturbations of Toda include the usual
$\alpha+\beta$ model. In this paper we explore and compare two
families, or hierarchies, of FPU models, closer and closer to either
the linear or the Toda model, by computing numerically, for each
model, the maximal Lyapunov exponent $\chi$. We study the asymptotics
of $\chi$ for large $N$ (the number of particles) and small $\eps$
(the specific energy $E/N$), and find, for all models, asymptotic
power laws $\chi\simeq C\eps^a$, $C$ and $a$ depending on the
model. The asymptotics turns out to be, in general, rather slow, and
producing accurate results requires a great computational effort. We
also revisit and extend the analytic computation of $\chi$ introduced
by Casetti, Livi and Pettini, originally formulated for the
$\beta$-model. With great evidence the theory extends successfully to
all models of the linear hierarchy, but not to models close to Toda.

\end{abstract}
\noindent
\vfill\eject

\section{Introduction}

\subsection{The purpose and the considered models}

This paper is devoted to the Fermi-Pasta-Ulam (FPU) problem
\cite{FPU}, precisely to a numerical study of the maximal Lyapunov
exponent $\chi$ in different one-dimensional FPU models. Numerical results will
be also critically compared with the existing theory. The Hamiltonian is
\beq{H} H(p,q)=\frac12\sum_{i=1}^{N-1}
p_i^2+\sum_{i=1}^{N}V(q_{i}-q_{i-1})\ , \qquad q_0=q_{N}=0\ ,
\eeq
$V$ being some nearest-neighbours potential with single minimum in zero.

Two integrable models of this form are known, namely the linear
model with purely quadratic potential $V_L(r)=\frac12r^2$, which
trivially separates into independent normal modes, and the Toda
model \cite{Toda} \cite{henon-toda} \cite{flaschka-toda} with potential
\beq{V-T}
V_T(r)=\ctoda^{-2}(e^{\ctoda r}-1-\ctoda r)=V_L(r)+\cO(\ctoda
r^3)\ ;
\eeq
no other FPU-like integrable systems do exist, see \cite{D08}.

As in most papers on the subject, we are mainly interested in the
dynamics at small energy per particle $\eps=E/N$, paying however as
much attention as possible to the asymptoticity in $N$ at fixed $\eps$
(that is, ideally, to the thermodynamic limit).

The FPU models we are going to consider include the two
most studied ones, namely:

\begin{itemize}
\item[i. ] The so-called ``$\alpha$+$\beta$'' model, with
  potential
\beq{V-ab}
  V_{\alpha+\beta}(r)=\frac{r^2}{2}+\alpha\frac{r^3}{3}+\beta\frac{r^4}{4}\ ,
  \qquad \alpha\ne0,\quad\beta>0\ .
\eeq
The model is generally thought of
as a perturbation of the linear model, the cubic and quartic term
providing a small coupling between normal modes. The nonlinearity (the
mode coupling) is given by the two quantities
$$
|\alpha|\sqrt{\eps}\ ,\qquad \beta\eps\ ;
$$
at small $\eps$ the former dominates and correspondingly the distance
(in a rough sense) to the linear model is
$$
\dist(\alpha+\beta,\mathrm{Lin})\sim|\alpha|\sqrt{\eps}\ .
$$
The \ab\ model, however, is tangent up to the order three to the Toda
model with $\ctoda=2\alpha$, the distance being, roughly,
$$
\dist(\alpha+\beta,\mathrm{Toda})\sim|\beta-\beta_T|\eps\ ,
\qquad \beta_T=\frac23\alpha^2 \ .
$$
The relevance of the underlying Toda dynamics in the \ab\  FPU model
was observed and emphasized already in 1982 by Ferguson, Flaschka and
McLaughlin in ref. \cite{FFM82} (a paper that we consider among the
fundamental ones, after the original FPU paper).

\item[ii. ] The symmetric ``pure $\beta$'' model, with potential
\beq{V-b}
V_{\beta}(r)=\frac{r^2}{2}+\beta\frac{r^4}{4}\ ,
\qquad \beta>0\ .
\eeq
Its closest integrable approximation is the linear
model and, roughly,
$$
\dist(\mathrm{pure}\ \beta,\mathrm{Lin})\sim \beta\eps\ .
$$
\end{itemize}

\noindent
Models with a higher order contact with Toda have been taken into
consideration, as far as we know, in \cite{CPC00} and in a more recent
paper by two of us \cite{BP11-1}. Such models are:

\begin{itemize}
\item[iii. ] The \ab\ model, with however
$\beta=\beta_T=\frac23\alpha^2$. We shall call it the ``\bt'' model.
Quite trivially, 
$$
\dist(\beta_T,\mathrm{Toda})\sim |\gamma_T|\eps^{3/2}\ ,\qquad
    \gamma_T=\frac13\alpha^3\ .
$$
    
\item[iv. ] The ``$\gamma_T$'' model, with 
$$
  V_{\gamma_T}(r)=\frac{r^2}{2}+\alpha\frac{r^3}{3}+\beta_T\frac{r^4}{4}
    +\gamma_T\frac{r^5}5+\delta\frac{r^6}6\ ,
  \qquad 0<\delta\ne\delta_T=\frac2{15}\alpha^4
$$
(the term of degree 6, with positive coefficient, ensures stability).
Roughly,
$$
\dist(\gamma_T,\mathrm{Toda})\sim|\delta-\delta_T|\eps^2\ .
$$

\end{itemize}

\noindent
A corresponding hierarchy of models closer and closer to the linear model
can obviously be produced; for completeness, and to test the existing
theory, we investigated two of them, namely

\begin{itemize}
\item[v. ] The ``\gd'' model, with potential
  \beq{V-gd}
  V_{\gamma+\delta}(r)=\frac{r^2}{2}+\gamma\frac{r^5}{5}+\delta\frac{r^6}{6}\ ,
  \qquad \alpha\ne0,\quad\delta>0\ ;
\eeq
it is obviously
$$
\dist(\gamma+\delta,\mathrm{Lin})\sim \gamma\eps^{3/2}\ .
$$

\item[vi. ] The ``\pd'' model, with potential
  \beq{V-pd}
  V_{\delta}(r)=\frac{r^2}{2}+\delta\frac{r^6}{6}\ ,
  \qquad \delta>0\ ,
\eeq
and correspondingly
$$
\dist(\mathrm{pure}\ \delta,\mathrm{Lin})\sim \delta\eps^{2}\ .
$$
\end{itemize}

\noindent
What is lacking in such a landscape of nearly integrable models is a model
distant both from the linear and from the Toda model by cubic terms in
$V$. The one we suggest here is
\begin{itemize}
  \item[vii. ] A generalization of \ab, with site-dependent potential
$$
V_i(r)=\frac{r^2}2+\alpha_i\frac{r^3}3+\beta\frac{r^4}4
$$
(see later for the choice of $\alpha_i$). We shall call it the
``variable $\alpha$'' model. The distances are
$$
\dist(\mathrm{var.}\,\alpha,{\mathrm{Lin}})\sim
\dist(\mathrm{var.}\,\alpha,{\mathrm{Toda}})\sim
\alpha\sqrt{\eps}\ ,
$$
$\alpha$ denoting here, roughly, the size of the $\alpha_i$'s and of
their oscillation. 
\end{itemize}

\noindent
Models i-vii form a kind of constellation around the two integrable
models, which is symbolically represented in figure \ref{F-1}. The aim
of this paper is to explore such a constellation by means of the most
flexible tool, sensitive to the lack of integrability and independent
of the  choice of special coordinates such as the normal modes (as is
especially important around Toda), namely
the computation of the maximal Lyapunov exponent. The paper, in a
sense, is a continuation of \cite{BP11-1} and \cite{BCP13} (which in
turn, in their aim, are continuations of \cite{FFM82}).

\begin{figure}
\centerline{\epsfig{figure=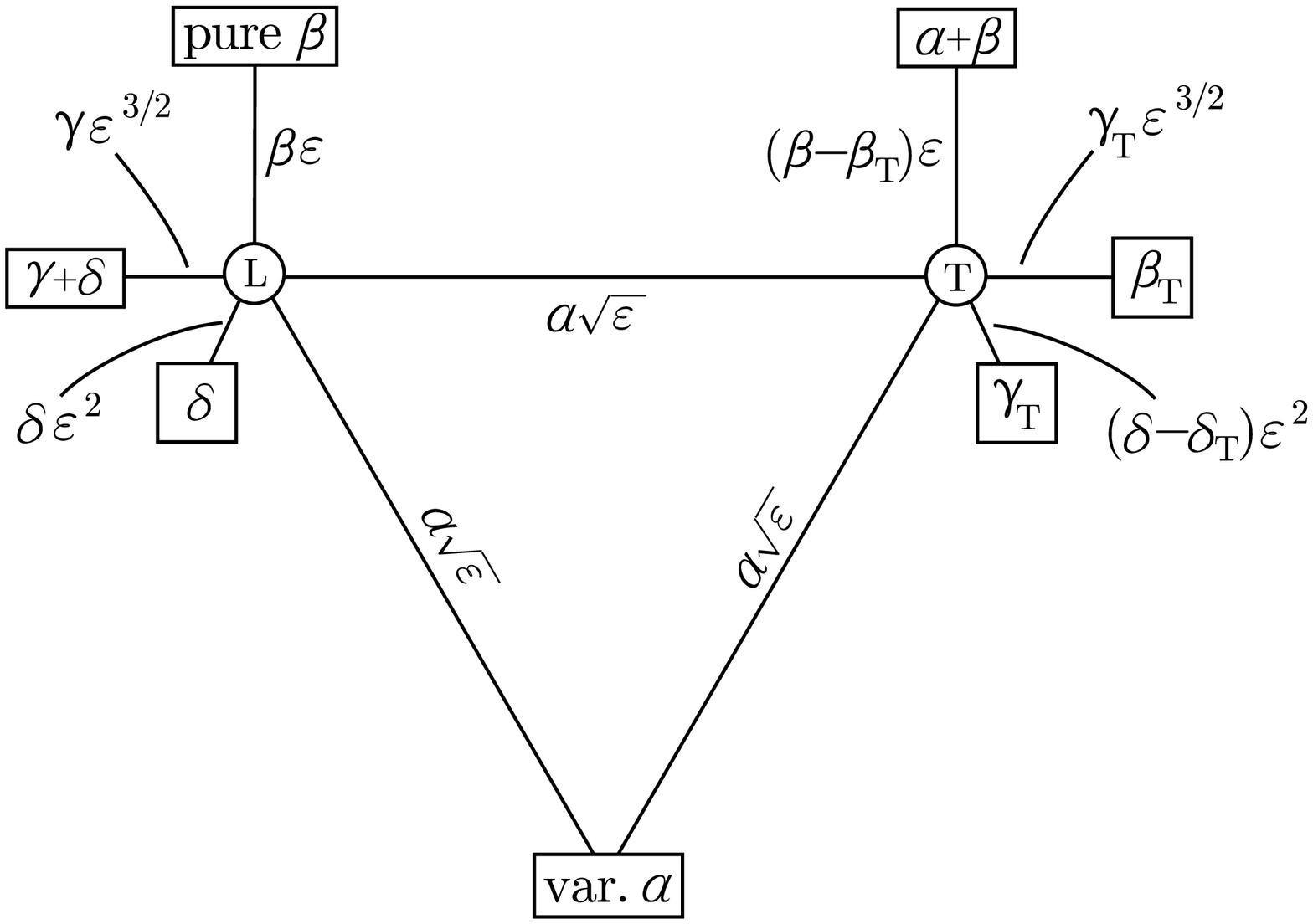,width=85mm}}
\Caption{F-1}{The different FPU models studied in this paper and their
rough distances to the linear (L) and to the Toda (T) models.}
\end{figure}

\bl\noindent {\it Remark.} For the \ab, \bt\ and \gt\ models, the
value of $\alpha$, if different from zero, is irrelevant and
simply fixes the energy scale: models with the same ratios
$\beta/\alpha^2$, $\gamma/\alpha^3$, $\delta/\alpha^4$... have
indeed identical dynamics up to trivial rescalings. We always used
$\alpha=-1$; correspondingly, $\beta_T=2/3$, $\gamma_T=-1/3$,
$\delta_T=2/15$. (Negative $\alpha$ means the ``springs'' between
masses get stiff when compressed, as in typical molecular potentials.)
Similar considerations hold for the leading nonlinearity constant of
any of the above models. 

\subsection{Main connections with the literature}

\noindent
{\it A. General references on FPU.} \quad The literature on FPU, over more
than 60 years, is so huge, that it is not even conceivable to
summarize it here. The state of the art, updated to 2005--06, can be
found in two collections of papers \cite{chaos_issue}
\cite{springer_fpu}, but it is hard---and this is a good indication of
the difficulty of the problem---to extract, from all of these papers
or partial reviews, a unitary view, free of contradictions.  In such
an uncommon situation, we can only indicate which are the main
references that have been important for us, to form our personal view
of the FPU problem, and should be considered at the basis of the
present work. (We are aware that the list is incomplete and 
important papers are possibly missing.)

After the original FPU paper, we wish to mention \cite{TM72}, where
the presence of longer time recurrencies was established, and some
early papers like \cite{IC66} \cite{IKC68} \cite{BSBL70}
\cite{GS72}, raising and discussing, with different perspectives, the very
crucial question about the persistence of the FPU paradox in the
thermodynamic limit.

A double turning point occurred in 1982: on the one hand, in the
already emphasized reference \cite{FFM82}, evidence was provided that
the FPU results (in particular, the partial energy sharing among low
modes in the \ab\ model) are fully explained by the underlying
integrable Toda dynamics; on the other hand, in two similarly crucial
papers \cite{Fucetal} \cite{LPRSV83}, the so-called two-time scales
(or metastability) scenario was introduced: in essence, on a first
relatively short time scale the phenomena observed by FPU (essential
lack of ergodicity and mixing, quasi-periodicity) occur, while on a
much larger time scale the normal statistical behavior is
recovered. The role of metastability was discussed and
emphasized in later papers like \cite{BGG04} \cite{BGP04} \cite{CGGP07}.

\begin{figure}
\centerline{\epsfig{figure=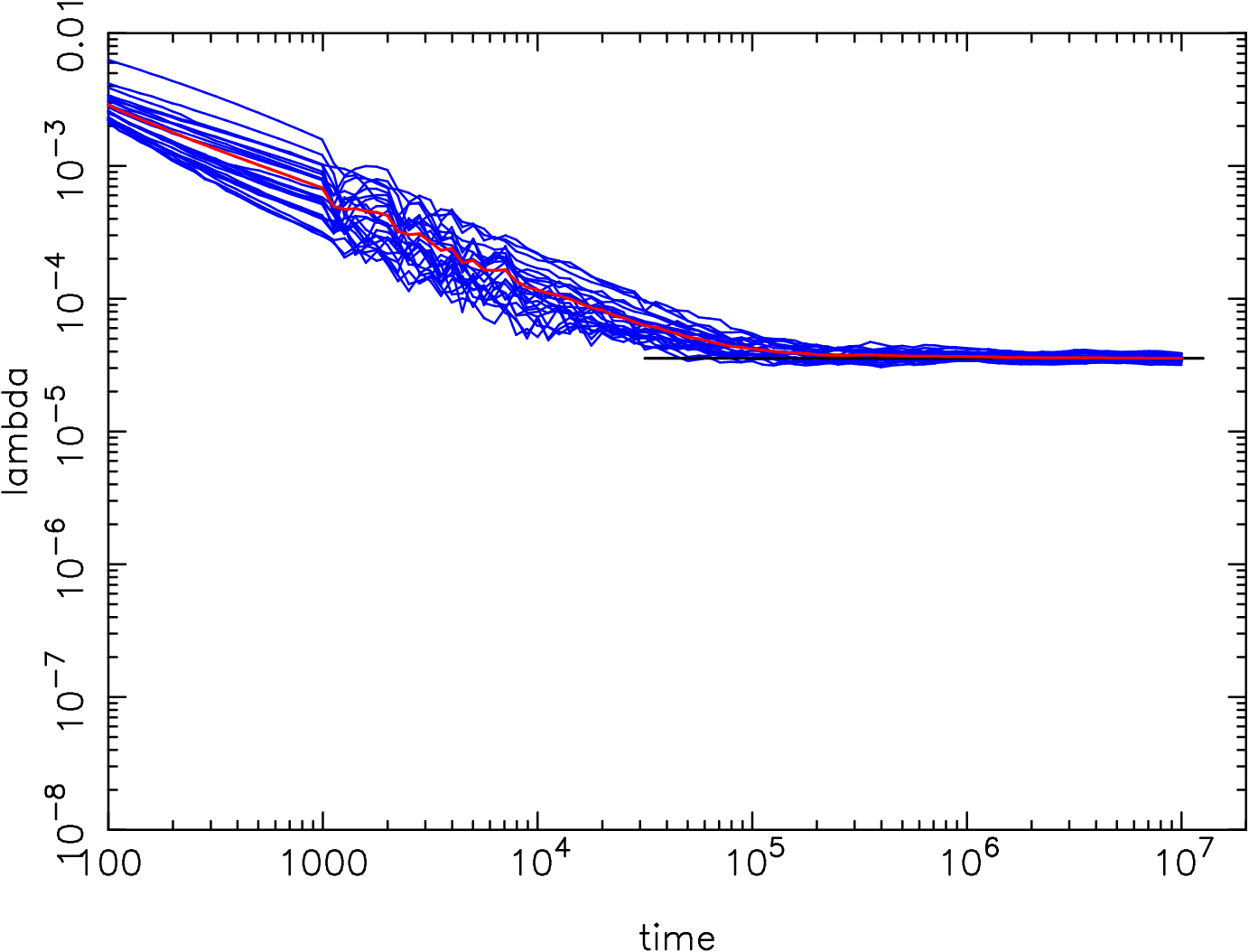,width=65mm,angle=-90}
\hspace{3mm}\epsfig{figure=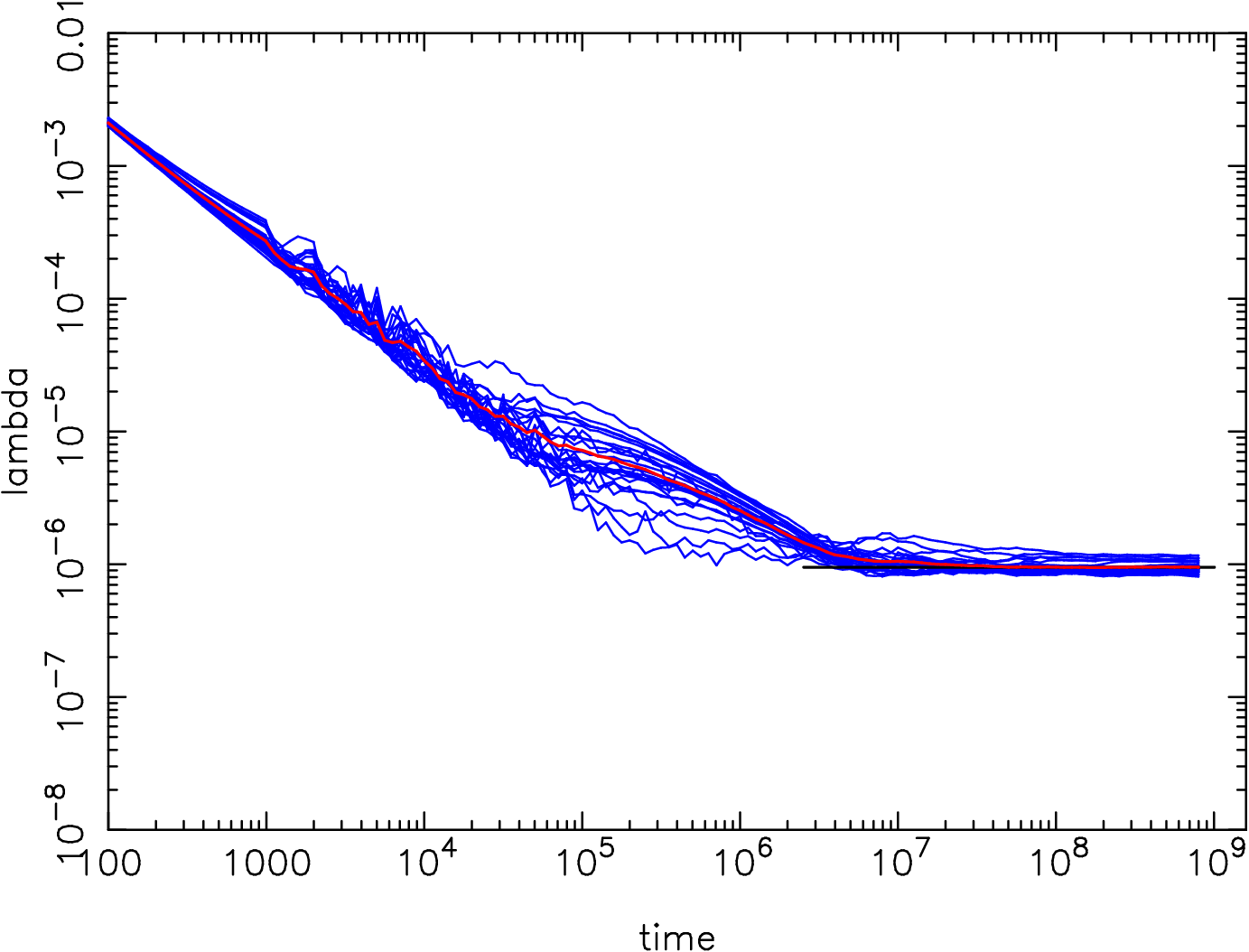,width=65mm,angle=-90}}
\vspace{3mm}
\centerline{\epsfig{figure=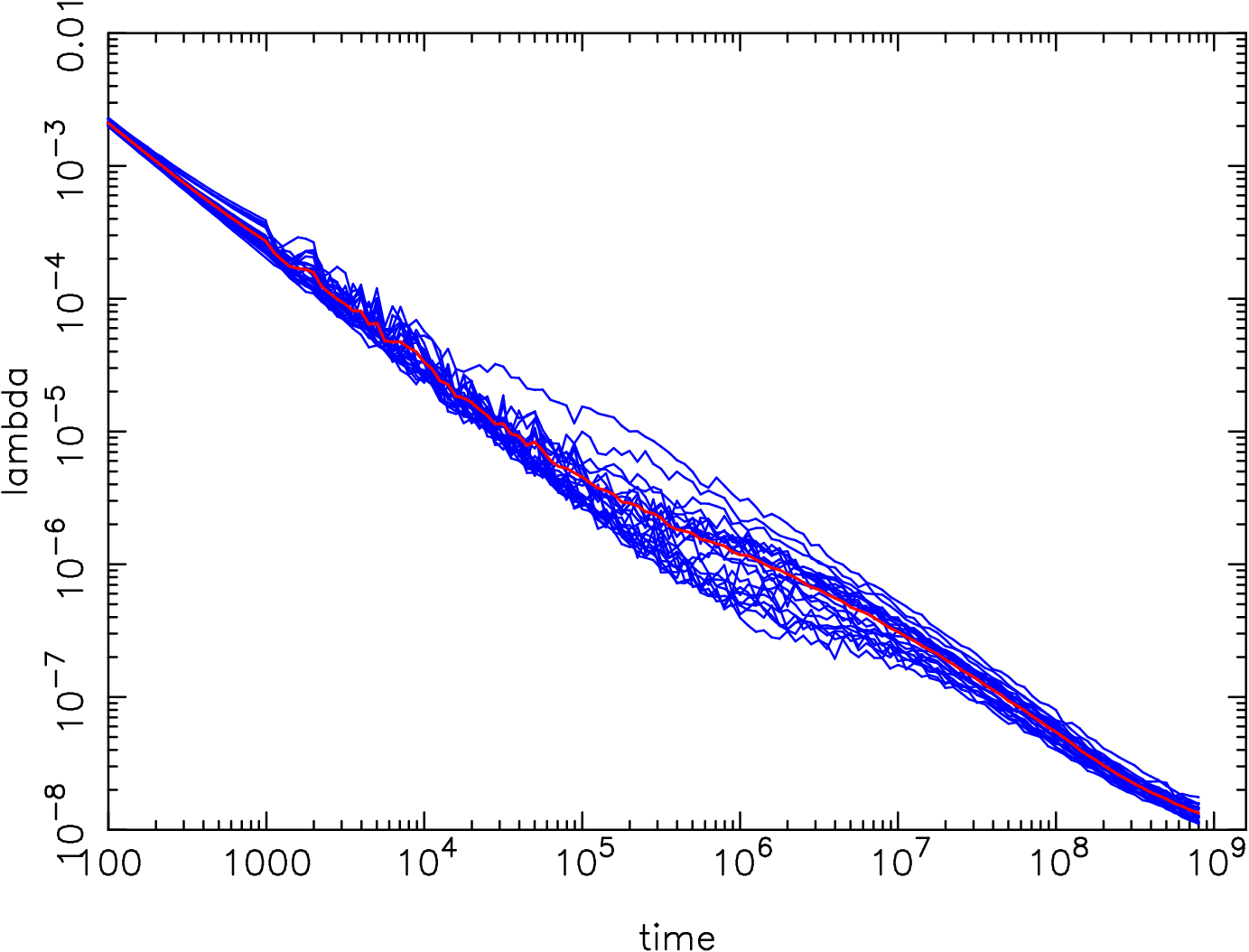,width=65mm,angle=-90}}
\Caption{F-2}{Upper panels: $\hchi(x,t)$ for the \ab\ model with $N=1024$,
  $\alpha=-1$ and $\beta=2$, at $\eps=10^{-3}$
  (left) and $10^{-4}$ (right); 24 random choices of the initial datum
  $x$ (blue) and their average $\chib(t)$ (red). The black horizontal
  segment is useful to appreciate the reaching of the limit ??). Lower
panel: same quantities for the Toda model, with $N=1024$, $\eps=10^{-4}$. }
\end{figure}

We expressed our view on the one-dimensional FPU problem in
\cite{BLP09}, \cite{BP11-1} and \cite{BCP13}. In short:

\begin{itemize}
\item[-- ] For models of the Toda hierarchy, i.e. \ab, \bt\ and \gt,
  the ideas introduced in \cite{FFM82} (the underlying integrable Toda
  dynamics) and in \cite{Fucetal} \cite{LPRSV83} (the two
  time--scales) naturally combine in a unitary quite elementary
  picture: in the shorter time scale of the metastability scenario the
  Toda actions stay almost constant, while the angles turn on tori
  producing a partial averaging (producing in particular the FPU state,
  if only a few low
  frequency modes are excited); on larger times the actions slowly
  diffuse throughout the phase space, with small speed for small
  $\eps$, producing eventually complete microcanonical averaging. The
  puzzling special behavior of low modes in FPU, during the short
  time scale, is due to a lack of correspondence between linear modes
  and Toda actions, when the index $k$ of the mode is, roughly,
  smaller than $\eps^{1/4}N$.

  \item[-- ] The thermodynamic limit, in all FPU models, is highly
    nontrivial: the limits $\eps\to0$ and $N\to\infty$ do not commute
    in an essential way, the asymptotics in $N$ becoming slower and
    slower at small $\eps$. Correspondingly, great attention must be paid
    in numerical computations.

  \item[-- ] In the thermodynamic limit, for small $\eps$, all time
    scales, including the large equilibrium ones, are in general
    (i.e., with the exception of special initial conditions
    \cite{BGP04} \cite{BP11-1}) power laws $\eps^{-a}$, the exponent
    $a$ depending on the degree of contact with the nearest integrable
    model, linear or Toda. In particular, with quite good accuracy,
    $a=9/4$ for the \ab\ and for the \pb\ models, $a=3$ for the
    \bt\ model and $a=4$ for the \gt\ model.
    
\end{itemize}

\noindent
{\it B. Specific references on the maximal Lyapunov exponent.} \quad
Let us now draw the attention to some keynote papers in which
Lyapunov Exponents are used to investigate the FPU problem. Here too,
the literature is rather abundant, and we limit ourselves to mention
the few papers which are more directly connected with the present
work. 

After the pioneering paper \cite{CDGS76}, the first explicit reference
to Lyapunov exponents in connection with FPU seems to be
\cite{BLT80}. Refs. \cite{PL90} and \cite{PC91} are devoted to the
\pb\ model (as well as to the so-called $\phi^4$ model), with $N$
between 32 and 128; a crossover between a weak and a strong
stochasticity regime is there made evident, and in the weak
stochasticity regime, at small $\eps$, the power law $\chi\sim \eps^2$
is found. Ref. \cite{CCPC97} contains instead a computation of $\chi$
in the pure $\alpha$ model (i.e. \ab\ with $\beta=0$), with $N$
between 8 and 64. In \cite{CPC00} the \bt\ model is studied again, with $N$
up to 1024; in the weak stochasticity regime the same power law
$\eps^2$ as in the \pb\ model is suggested.

A very remarkable sequence of papers, including in particular
\cite{CLP95} and \cite{CCP96} (see also sections 2.4.7 and 2.8.2 of the
review paper \cite{LLPR08}), develop a statistic--geometric theory, based on a
subtle analysis of the fluctuations of the curvature in the phase
space in a natural metrics, which leads to an analytic estimate of
$\chi$ as a function of $\eps$, in the thermodynamic limit. The
theoretical curve fits impressively well the numerical data of the
\pb\ model, both in the small $\eps$ and in the large $\eps$
regimes. Asymptotically for small $\eps$ the power law
$$
\chi = C\,\eps^2
$$
(with computed $C$) is found.

Further comments on the theory, including a natural extension of it
to all models in the hierarchy around the linear model (but not
around Toda, where the theory apparently does not apply) are demanded
to Section 3.

\subsection{A sketch of results}
Our results can be summarized by saying that:

\begin{itemize}
\item[-- ] For all models, at large $N$ and small $\eps$, the Lyapunov
  exponent obeys a power law $\chi=C\,\eps^a$, $C$ and $a$ depending
  on $N$; for large $N$, the exponent $a$ always converges (although
  possibly very slowly) to a limit value, while the constant $C$, for
  the most important models namely \ab\ and \pb, possibly mantains a
  residual logarithmic dependence on $N$.

\item[-- ] Models of the linear hierarchy exhibit ``simple''
  exponents, namely $a=2$ (in agreement with the above quoted
  literature) for the \pb\ model, $a=3$ and $4$ respectively for the
  \gd\ and the \pd\ model, and $a=1$ for those \av\ models which, see
  Section 2.4, are close to the linear model.

\item[-- ] Models of the Toda hierarchy exhibit less understandable
  exponents, actually $a=3/2$, $1.9$ and $2.2$ respectively for the \ab,
  the \bt\ and the \gt\ model; $a=1.2$ for some other \av\ models, which are
  close to Toda.

\item[-- ] The theory developed in \cite{CLP95} \cite{CCP96} in
  connection with the \pb\ model successfully extends to all models of
  the linear hierarchy, producing, for large $N$ and small $\eps$, the
  above mentioned values of $a$. The theory instead, as far as we
  understand it, does not extend to the models of the Toda hierarchy
  and to Toda itself (for which of course $\chi=0$). The obstruction
  seems precisely to be the difficulty to detect the integrability of
  Toda, and consequently the peculiarities of the nearby models.

\end{itemize}

\noindent{\it Akcnowledgements}: we are indebted to Roberto Livi
(Firenze) and to Luigi Galgani and Andrea Carati (Milano) for helpful
stimulating discussions. 

\section{Numerical results}

\subsection{Computing $\chi$}

Trajectories have been computed by a symplectic integration algorithm,
namely a fourth-order leap-frog based on 
\cite{Y90}. For general questions concerning symplectic integration
algorithms see \cite{BG94}, \cite{H94}, \cite{HLW02}; for their use in
connection with FPU models see \cite{PP05}, \cite{BP11-2}.
The Lyapunov exponent $\chi$ has been computed by the traditional
algorithm \cite{BGS} \cite{BGGS80}. For any $x=(p,q)$ and any tangent
vector $\xi=(\delta p,\delta q)$, let $\Phi^t$ and $D\Phi_x^t$ denote
the flow and, respectively, its tangent application in $x$, and let
\beq{chi-xt}
\hchi(x,t)=\frac1t \log \frac{\|D\Phi_x^t\xi\|}{\|\xi\|}
\eeq
be the usual finite-time Lyapunov exponent. The limit for $t\to\infty$ of
$\hchi(x,t)$ gives, for (almost) any choice of $\xi$, the maximal
Lyapunov exponent in $x$. Practically, the dependence on $\xi$ at the
r.h.s. of \calleq{chi-xt} is soon
lost, much before the limit, and for this reason $\xi$ does not appear
among the arguments of $\hchi$. The choice of the norm in the
tangent space is similarly irrelevant (we used the Euclidian norm in
the $\delta p,\delta q$ coordinates).

For any model, any $N$ and any $\eps$, we found convenient to consider
24 different initial data $x$, extracted randomly on the chosen
constant energy surface with a Gaussian measure on normal modes (a
measure close to the microcanonical one).  As an example, figure
\ref{F-2}, upper panels, shows $\hchi(x,t)$ as a function of $t$, as
usual in log-log scale, for the \ab\ model with $N=1024$, $\beta=2$,
at $\eps=10^{-3}$ (left) and $\eps=10^{-4}$ (right); the blue curves
are the 24 individual curves $\hchi(x,t)$, the red curve is their
average $\chib(t)$ (the average on the 24 initial data). Quite
clearly, the finite time behavior depends on $x$ (due to the different
paths in the phase space), but the limit is the same. Correspondingly,
the precise choice of the 24 random initial data is quite irrelevant,
and averaging is only helpful to produce a cleaner curve. Practically,
different samples of the 24 initial data produce hardly
distinguishable averages $\chib(t)$ much before the limit.  For all
models and any $N$ and $\eps$, the behavior of $\hchi(x,t)$ and
$\chib(t)$ is as in figure \ref{F-2}. The {\it plateau} is obviously
absent in Toda, for which (if the numerical integration is accurate,
see Section 4.A) the maximal Lyapunov exponent vanishes and
correspondingly both $\hchi(t,x)$ and the average $\chib(t)$ decrease
indefinitely, see the lower panel of figure \ref{F-2}.

We shall denote by $\chi(N,\eps)$ the
limit value of $\chib(t)$. This is the quantity we are going to investigate.

\begin{figure}
\centerline{\epsfig{figure=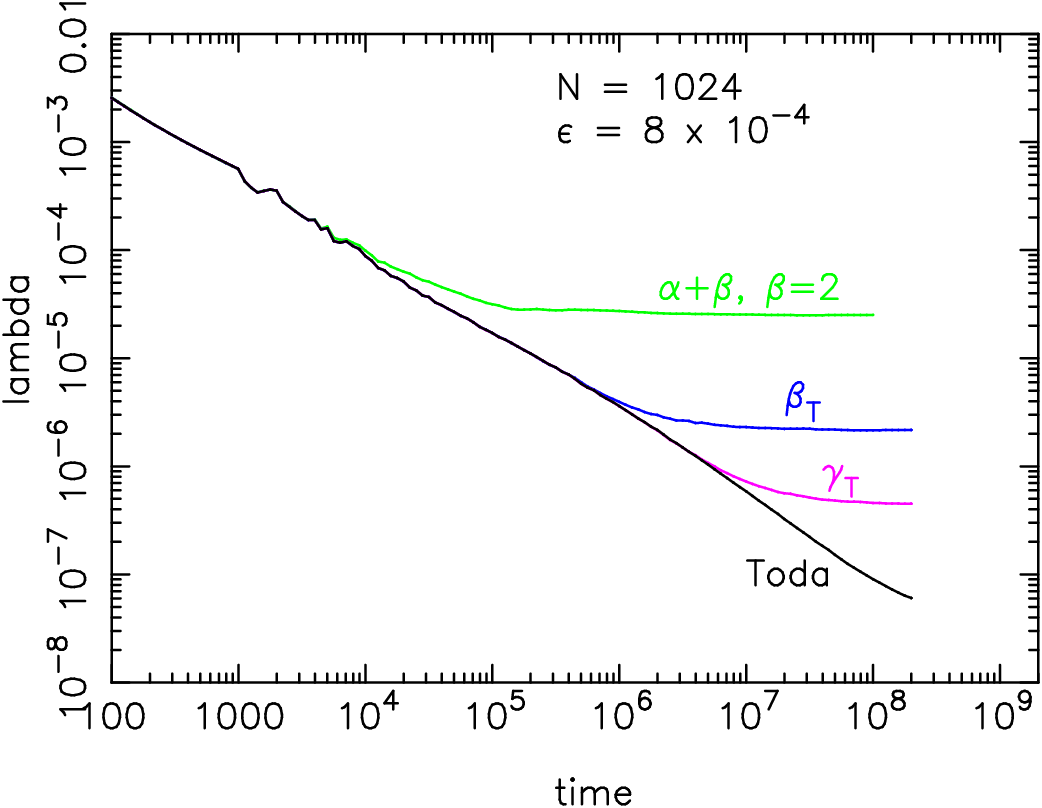,width=90mm} }
\Caption{F-3}{$\chib(t)$ as function of $t$ for the models of the Toda
  hierarchy. Models, top to bottom: \ab\ with
  $\beta=2$, \bt, \gt, Toda. Same $N=1024$, $\eps=8\times10^{-4}$.}
\end{figure}

\subsection{$\chi$ for the Toda hierarchy}

\noindent
{\it A. Generalities.} \quad
In this subsection we address to models close to Toda, namely \ab,
\bt\ and \gt.
A preliminary indication that the proximity to Toda is relevant is
provided by fig. \ref{F-3}, where $\chib(t)$ is plotted vs. $t$ for
these models as well as for Toda (same initial data). As is not surprising,
the closer a model is to Toda, the longer its $\chib(t)$ follows
$\chib(t)$ of Toda.

\noindent
{\em Remark.} As a general well known fact, the finite time Lyapunov exponent
$\hchi(x,t)$, in any nearly-integrable system, follows closely the
corresponding curve of the reference integrable model, until it detaches
and reaches a {\it plateau} (sometimes from above as in figure \ref{F-3},
sometimes from below). Using the average $\chib(t)$ in place of
$\hchi(x,t)$ cleans the curves and makes the phenomenon more
evident (look how the curves in figure \ref{F-3} superimpose to the
Toda curve, even in the details). For FPU, this phenomenon is
commented for example in \cite{CCPC97} and, with main attention to the case
of small
$N$, in \cite{GPP05}. We propose here a quite elementary interpretation
(making reference, in the notation, to our $\chib(t)$).
In an integrable system, the norm of the tangent vector,
supposing it be one at $t=0$, grows roughly linearly in time, say
$n(t)\simeq 1+ht$, $h$ being some norm of the Hessian of the
Hamiltonian in action-angle variables. For a nearly-integrable system
it is, again roughly, $n(t)\simeq 1+ht+c(e^{\chi t}-1)$,
$\chi$ small, and correspondingly
\beq{cross}
\chib(t)\simeq\frac1t \log[1+ht+c(e^{\chi t}-1)]\ .
\eeq
Quite clearly, for $t<\chi^{-1}$ the behavior is almost identical to
the corresponding integrable system, and $\chib(t)$ decreases in time
as $\log(1+ht)/t$. At $t\simeq\chi^{-1}$ there is a crossover to a
{\it plateau} and asymptotically
$$
\chib(t)\simeq\chi+t^{-1}\log c \ ,
$$
the {\it plateau} being reached from above for $c>1$, from below for
$c<1$. $\chib$ as in \calleq{cross} is plotted vs. $t$ in figure
\ref{F-4}; the upper curve has $c>1$, the lower one has $c<1$. The
profiles are very typical of the computation of $\chi$, in
nearly-integrable models. Let us remark that our interpretation does
not need assuming that trajectories move, after a certain ``trapping
time'', from a regular to a chaotic region. We do not claim our view
based on \calleq{cross} is original, but we are not aware of a paper
where it is suggested. Let us further comment that, in this
perspective, the detachment time is not an
independent information on the dynamics, such a time essentially
coinciding with $\chi^{-1}$.

\begin{figure}
\centerline{\epsfig{figure=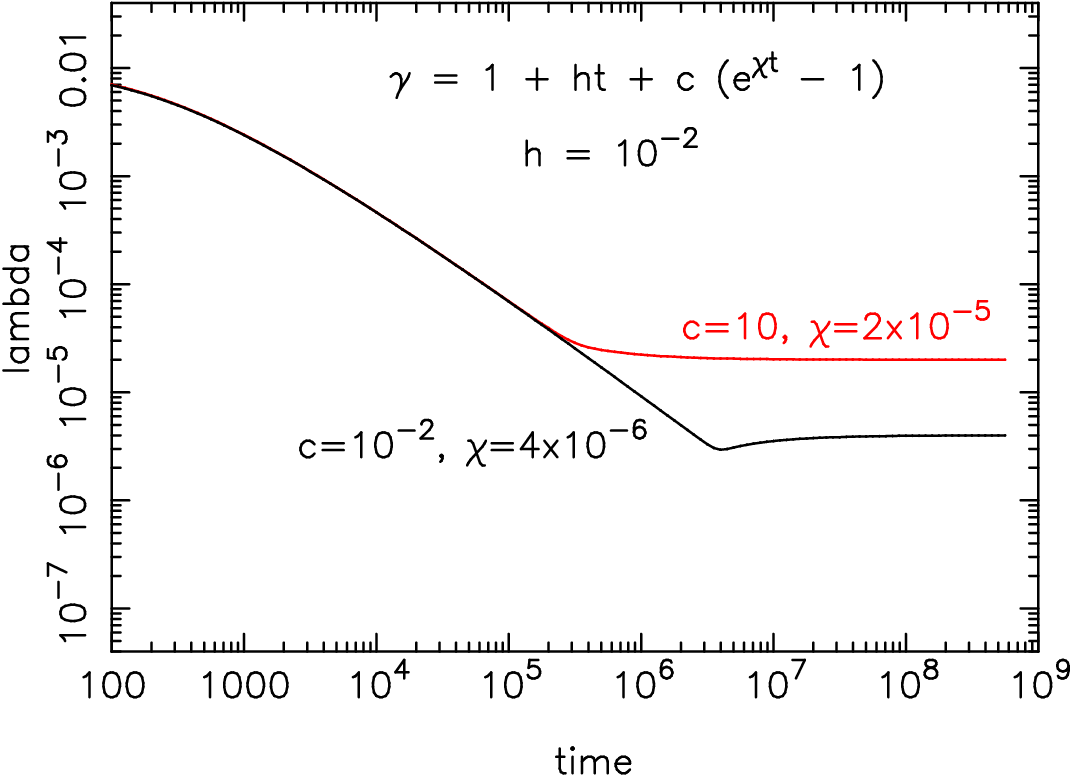,width=90mm} }
\Caption{F-4}{A plot of \calleq{cross}, for $h=10^{-2}$ and
  $c=10$, $\chi=2\times 10^{-5}$ (upper curve), $c=10^{-2}$,
  $\chi=4\times 10^{-6}$ (lower curve).}
\end{figure}

\begin{figure}
\epsfig{figure=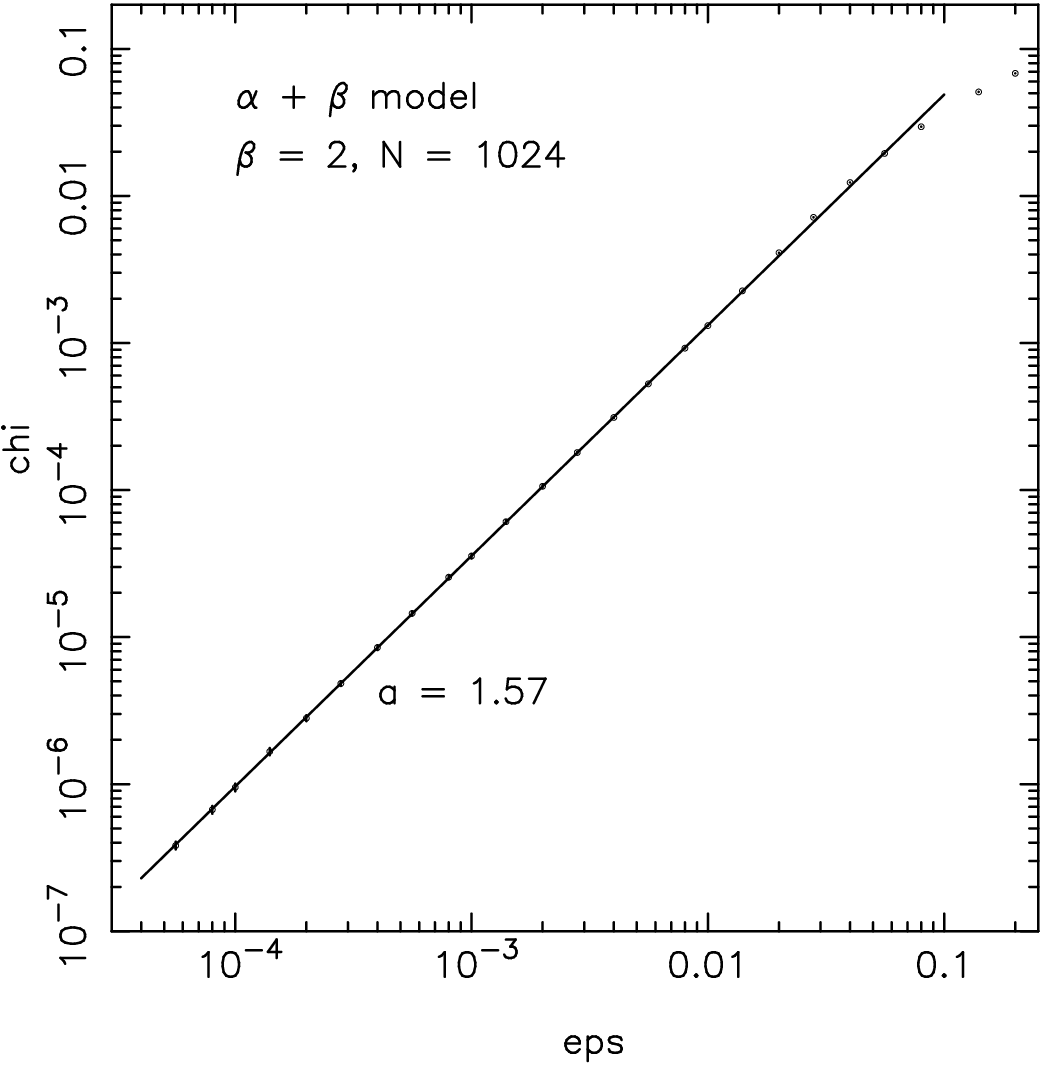,width=80mm}
\hspace{3mm}
\epsfig{figure=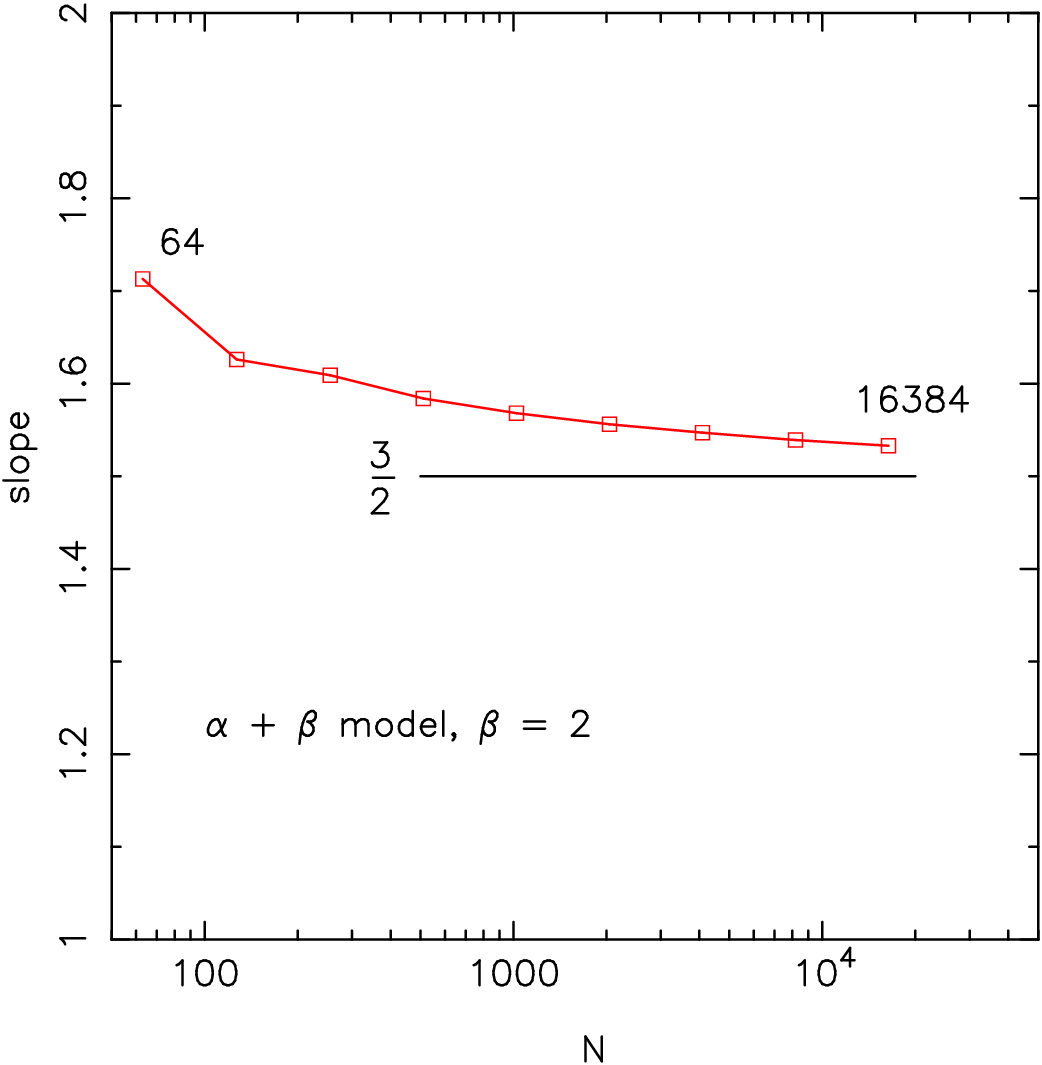,width=80mm}
\vspace{1mm}
\centerline{\epsfig{figure=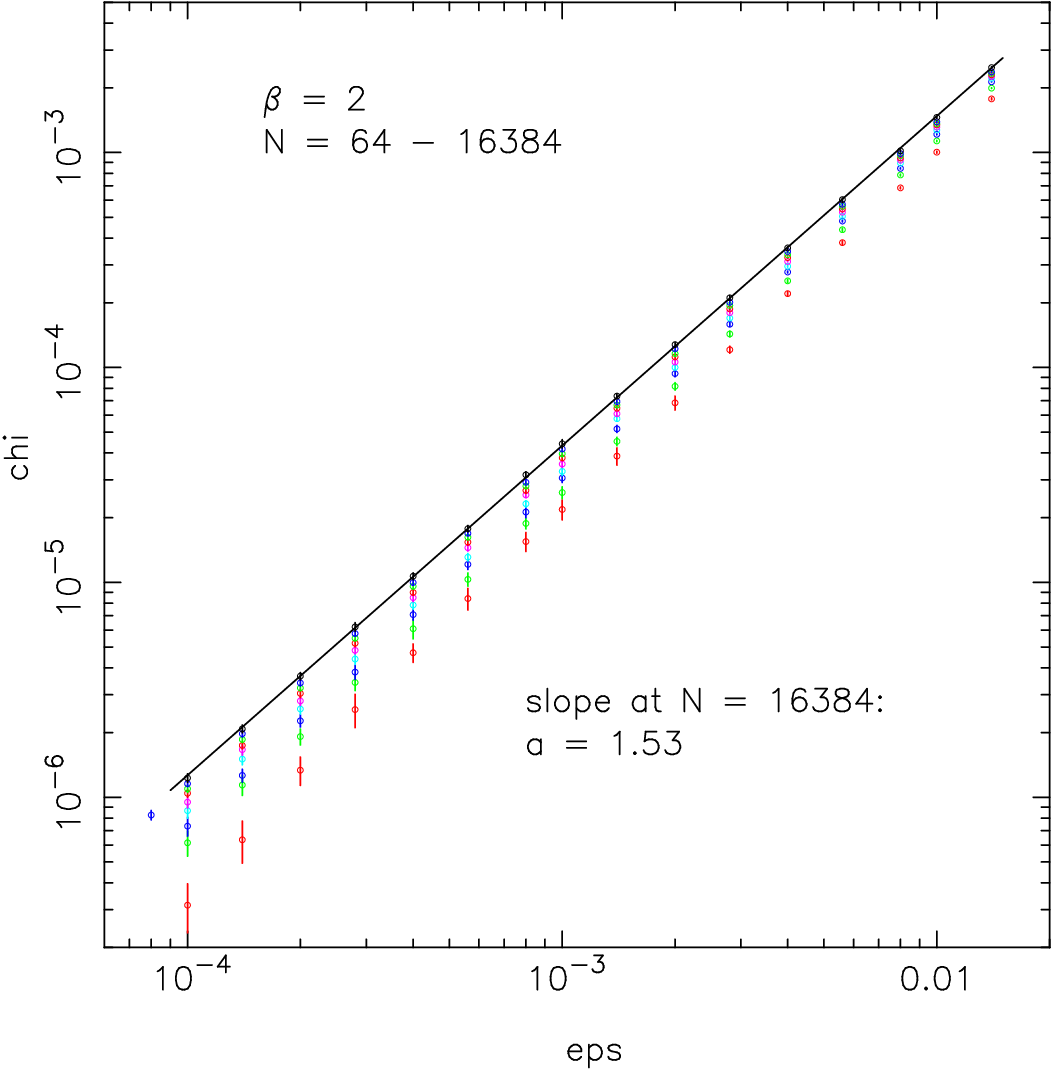,width=80mm}}
\Caption{F-5}{Up left: $\chi(N,\eps)$ as function of $\eps$ for the \ab\ model
  with $\beta=2$ and $N=1024$. Up right: the exponent $a$ as function of
  $N$, at fixed $\beta=2$ (least squares for
  $\eps<2\times10^{-2}$). Down: $\chi(N,\eps)$ as function of $\eps$,
  for $N=64,128, \ldots, 16\,384$. The line interpolates data 
  for $N=16\,384$.}
\end{figure}

\begin{figure}
\epsfig{figure=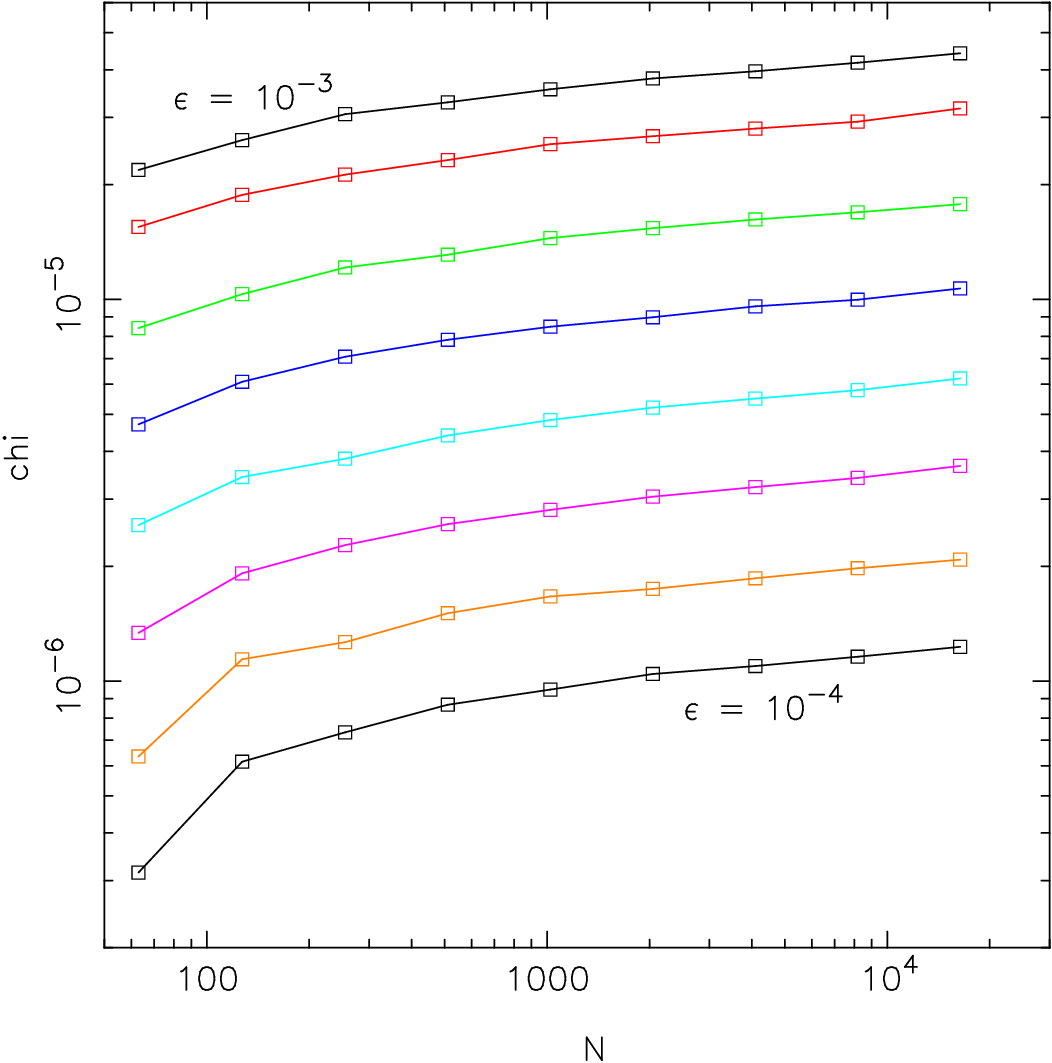,width=80mm}
\hspace{3mm}
\epsfig{figure=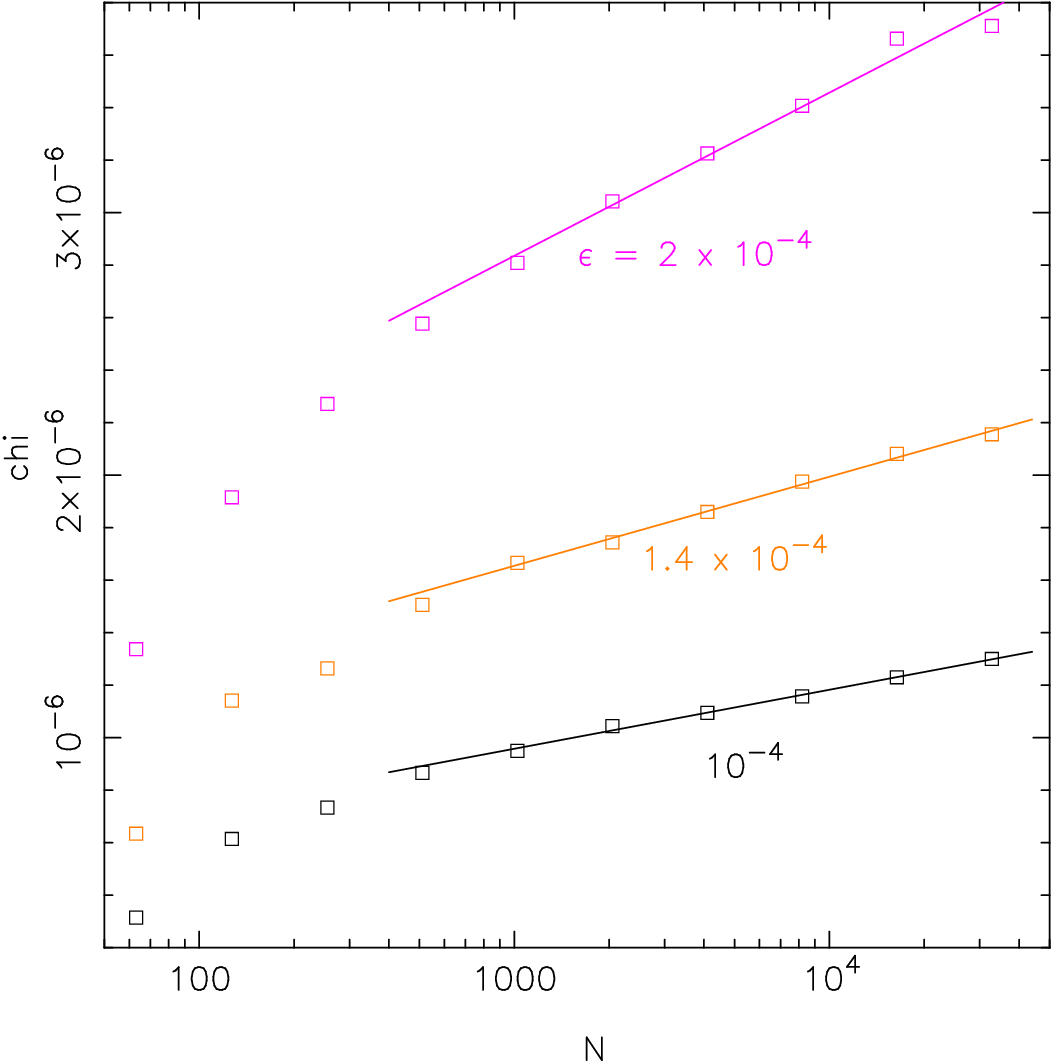,width=80mm}
\Caption{F-6}{$\chi(N,\eps)$ as function of $N$, for some values
  of $\eps$, in the \ab\ model with $\beta=2$. Left, log--log scale;
  right, selected values in semilog scale.}
\end{figure}

\bl\noindent {\it B. The \ab\ model.}\quad
Figure \ref{F-5}, upper right panel, 
shows $\chi(N,\eps)$ as function of $\eps$ for the
\ab\ model with $\beta=2$, $N=1024$. The data for $\eps<2\times
10^{-2}$ are quite well fitted by a line,\footnote{In principle, the
  data in the figure include an error bar, deduced in a standard way
  from the dispersion of the 24 individual data: let
  $\chi_1, \ldots, \chi_n$, $n=24$, be the values of $\hchi(x,t)$ for
  the $n$ considered trajectories at the largest computed $t$; then the
  error bar corresponds to $\pm3\sigma/\sqrt{n-1}$, $\sigma$ being the
  standard deviation of the $\chi_i$'s. Bars are practically invisible
  because of the same size of the symbols or smaller. }  
suggesting, as in \cite{PL90} \cite{PC91} \cite{CPC00} \cite{CLP95},
a power law
\beq{power}
\chi\sim\eps^a\ .
\eeq
The value of the exponent (least-squares best fit) is $a=1.57$, but the value 
depends on $N$ and, as shown in the upper right panel of the
figure, for large $N$ it seems to approach $a=3/2$. The
asymptotics looks unespectedly slow. The lower panel provides an overview of
$\chi(N,\eps)$ for several values of $N$ between $64$ and $16\,384$; 
values are less accurate for smaller $N$, as is shown
by the larger now well visible error bars.

In fact, some criticism is here mandatory, for the very existence of a
limit $N\to\infty$ for $\chi(N,\eps)$ at fixed $\eps$ is questionable. Figure
\ref{F-6}, left panel, shows $\chi(N,\eps)$ as a function of $N$,
log-log scale, for $\eps$ between $10^{-3}$ and
$10^{-4}$. The curves, at least at the values of $N$ we reached, seem
to mantain a residual dependence on $N$ (in spite of the fact that the
slopes $a(N)$, as shown in fig. \ref{F-5}, apparently attain a limit).
The right panel of figure \ref{F-6} reports a few of such 
curves in semilog scale, in the attempt to put in
evidence a possible logarithmic dependence on $N$. The interpolation
with lines, for $N\ge 10^3$, looks acceptable. (Of course, a power law
with a small exponent also fits the data).
The miss of the thermodynamic limit by a residual
logarithmic dependence on $N$, in our opinion, is welcome. We are not
aware of a similar situation, at least not in connection with
FPU. Great computational effort has been necessary to produce data as
accurate as in figs. \ref{F-5} and \ref{F-6}. 

\begin{figure}
\centerline{\epsfig{figure=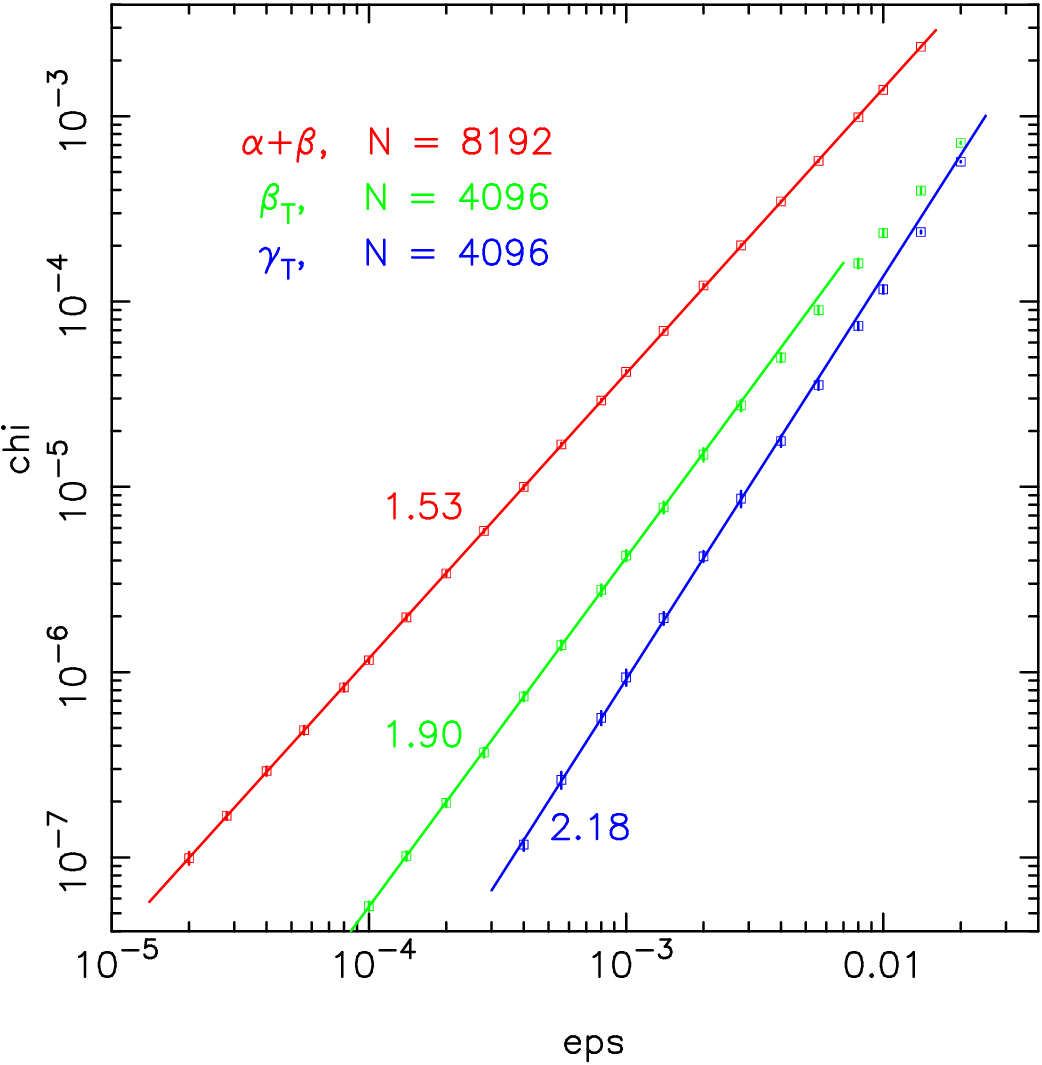,width=80mm}}
\Caption{F-7}{The different power laws for models of the Toda
  hierarchy ($N=8\,192$ for \ab, $n=4\,096$ for \bt\ and \gt).}
\end{figure}

\bl\noindent 
{\it C. Other models of the Toda hierarchy.}\qquad 
Computations have been repeated for the two mentioned models of the
Toda hierarchy, namely \bt\ and \gt, with $N$ up
to\footnote{Obtaining accurate results for such
  models is rather painful: the convergence of $\chib(t)$ to its
  limit is slow and requires very long computations. For this reason,
  $N$ was limited to $4\,096$. } 
$4\,096$. Qualitatively the behavior is similar to the
\ab\ model, $\chib$ obeying power laws like \calleq{power}.
 Figure \ref{F-7} summarizes the results, showing
together the different power laws we found for the different models.
In detail: 

\begin{itemize}
\item[$\circ$ ] The \bt\ model: we found $a\simeq 1.90$, a value close
  to $a=2$ approximately found, and conjectured as true, in
  \cite{CPC00}. We obviously agree that $a=2$ is an attracting
  conjecture, but such a value does not seem compatible with the
  numerical data. This is shown in figure \ref{F-8}, where
  $\chi(N,\eps)$ is reported as a function of $\eps$ for $N=256$,
  $512$, \ldots,  $4\,096$ (only the data for small $\eps$
  between $10^{-4}$ and $10^{-3}$ are reported). For each $N$, data
  look definitely well interpolated by lines, with slopes oscillating
  (without order) between $1.88$ and $1.90$. By comparison, a dashed
  line with slope 2, through the last point at $N=4\,096$, is also
  drawn. The slopes look well distinguishable. 

\item[$\circ$ ] The \gt\ model: we found $a\simeq 2.2$. Here too, the closest
  ``attracting'' value is $a=2$, but such a value looks not
  compatible with the numerical data. (In any case, it would be hard to
  conjecture $a=2$ for both the \bt\ and the \gt\ models.)

\end{itemize}

\noindent
The overall situation looks rather puzzling. We shall further discuss
this point among the Concluding remarks.

\begin{figure}
\centerline{\epsfig{figure=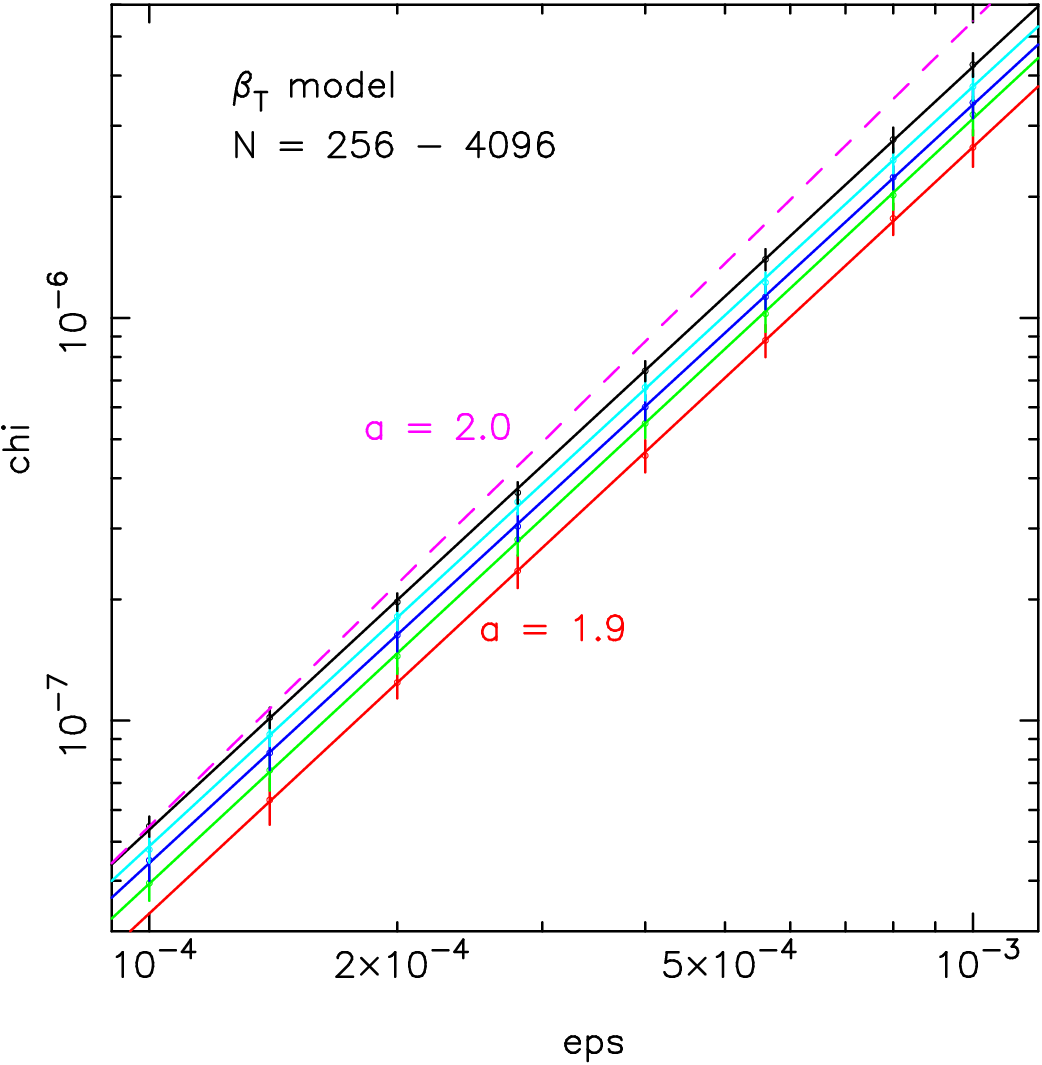,width=80mm}}
\Caption{F-8}{ Symbols and solid interpolating lines: $\chi(N,\eps)$
  as function of $\eps$, for (bottom to top)
  $N=256$, $512$, \ldots, $4\,096$. Dashed line: a line with slope $a=2$,
  through the last point at $N=4\,096$.}
\end{figure}

\subsection{$\chi$ for the hierarchy around the linear model}

\noindent
{\it A. The \pb\ model}\quad
We studied the power law \calleq{power} for the \pb\ model, although
established in the literature, so as to increase the accuracy of the
existing results. As already remarked in the Introduction, both the
available numerical data \cite{PL90} \cite{PC91} and the theory
\cite{CLP95} \cite{CCP96} suggest $a=2$. We confirm such a result,
asymptotically for large $N$, with however quite slow asymptotics in $N$, as
for the \ab\ model.

Figure \ref{F-9} shows, in the left panel, $\chi(N,\eps)$ as a
function of $\eps$, for $N=128$, $256$, \ldots, 
$8\,192$. The line interpolates data for $N=8\,192$. Its slope is
$a=2.05$, but as shown in the right panel of the figure, the slope
depends on $N$, and for $N\to\infty$ it likely converges to
$2$. Nevertheless, a residual dependence of $\chi$ on $N$,
approximately logarithmic as for the \ab\ model, apparently remains:
see figure \ref{F-10}, similar to figure \ref{F-6}. 

\begin{figure}
\epsfig{figure=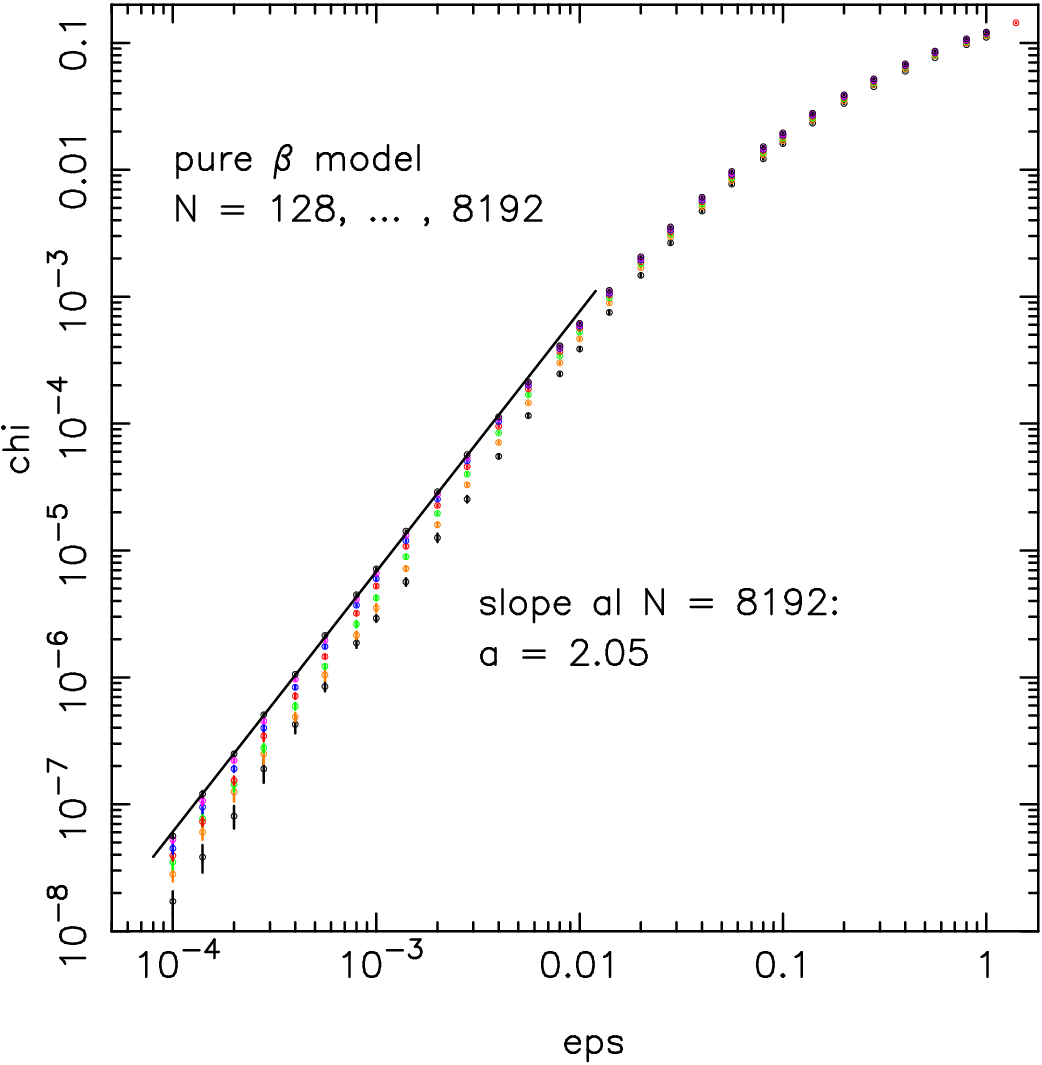,width=80mm}
\hspace{3mm}
\epsfig{figure=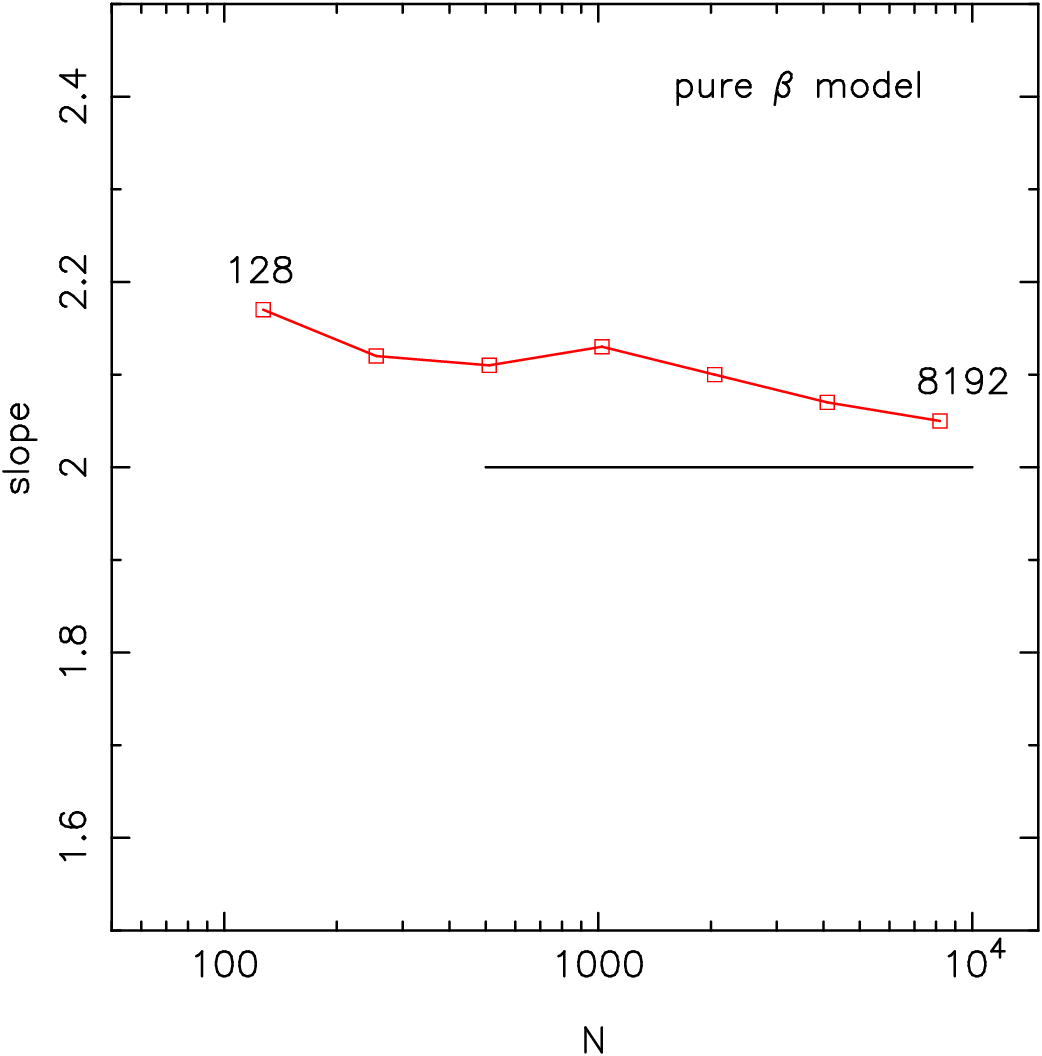,width=80mm}
\Caption{F-9}{
  Left: $\chi(N,\eps)$ as function of $\eps$, for the pure $\beta$ model,
  for $N=128$, $256$, \ldots, $8\,192$; 
  the line with slope $a=2.05$ interpolates data for $N=8\,192$. Right:
  the computed slope as function of $N$, suggesting a possible
  convergence to $a=2$.}
\end{figure}

\begin{figure}
\epsfig{figure=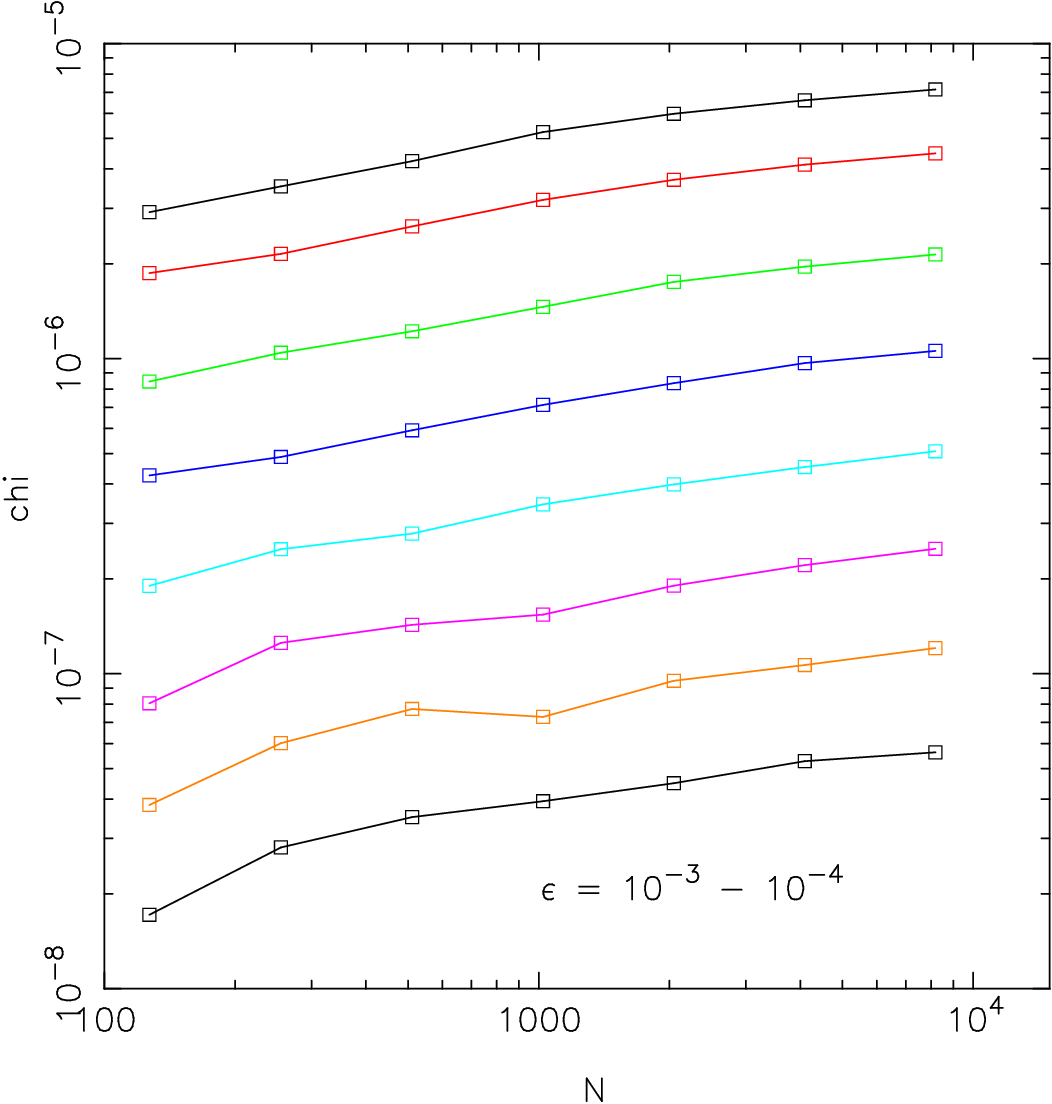,width=80mm}
\hspace{3mm}
\epsfig{figure=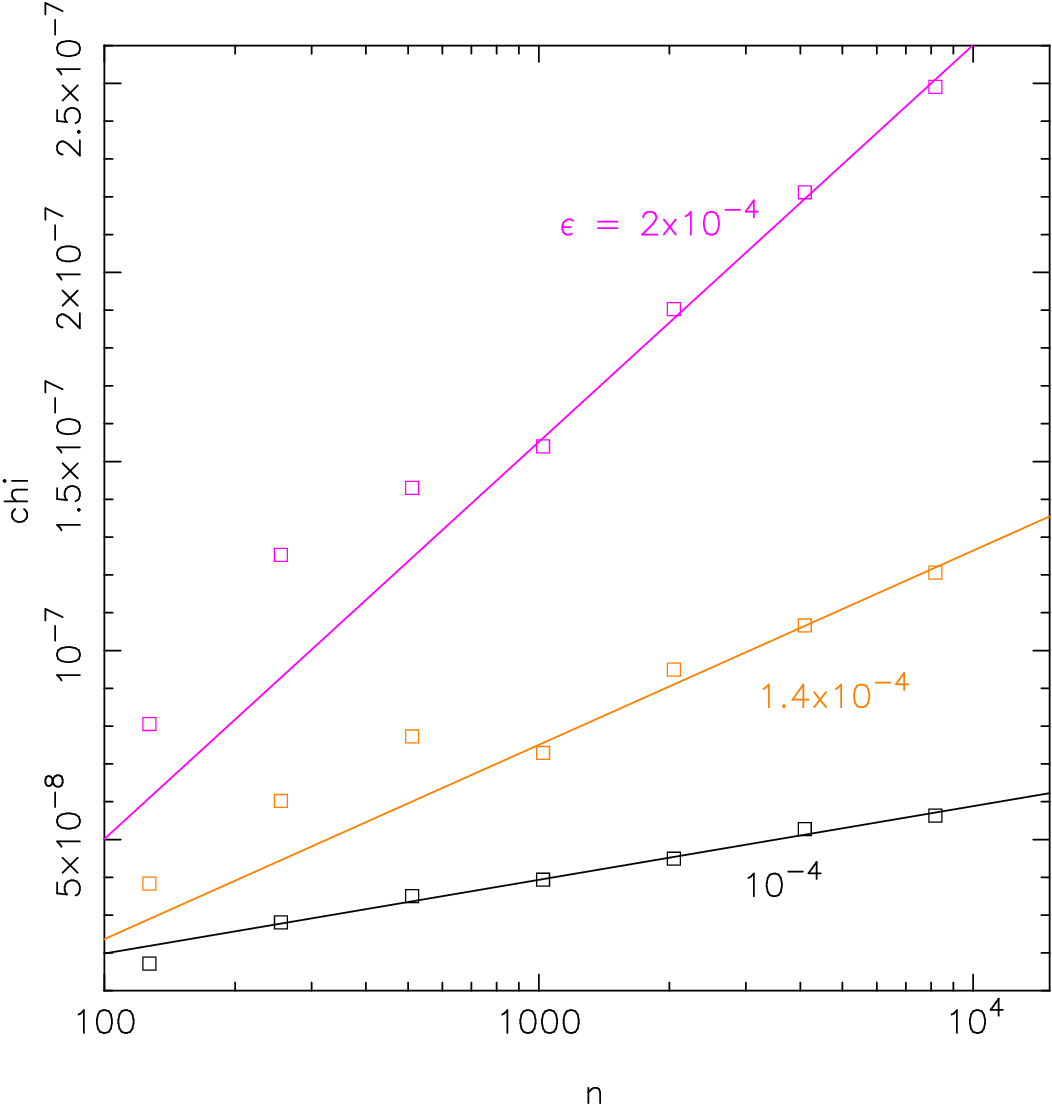,width=80mm}
\Caption{F-10}{$\chi(N,\eps)$ as function of $N$, for some values
  of $\eps$, in the pure $\beta$ model. Left, log--log scale;
  right, selected values in semilog scale.}
\end{figure}

\bl\noindent
{\it B. The next models.} \quad

As remarked in the Introduction, we investigated the two next models
of the hierarchy around the linear model. The main reason, besides the
pleasure of the symmetry of figure \ref{F-1}, is to make a test of the
theory developed in \cite{CLP95}, \cite{CCP96}.

Concerning the \gd\ model, the choice of the constants is $\gamma=1$
and $\delta=0.8$. The result is a
beautiful power law with $a=3$, essentially independent of $N$ at
least between $N=1024$ and $N=8\,192$; see figure \ref{F-11}, upper curve.
Concerning instead the \pd\ model, the natural choice is
$\delta=1$; for such a model we found $a=4$, see the lower curve of
the same figure. The very simple rule for the exponent, also including
the \pb\ model, seems to be
\beq{a-lh}
\chi\sim \eps^{s} \qquad \mathrm{for} \qquad
V(r)=\frac12 r^2+\cO(r^{2+s})\ ,
\eeq
for $s\ge 2$. As we shall comment in Section 3, such a rule is
in agreement with the theoretical prediction. 

\begin{figure}
\centerline{\epsfig{figure=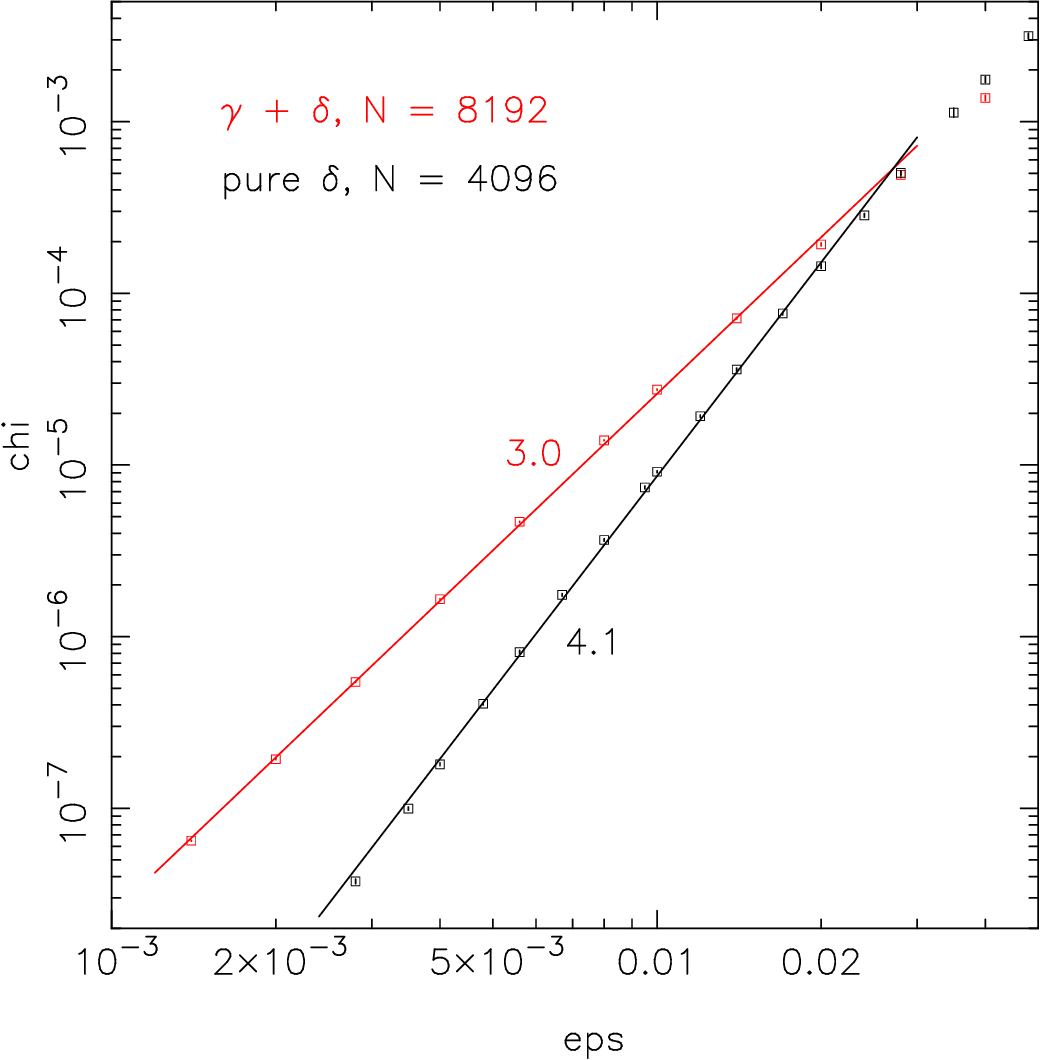,width=80mm}}
\Caption{F-11}{$\chi(N,\eps)$ as function of $\eps$ for the
  \gd\ model at $N=1024$ (upper curve), and the \pd\ model at $N=8192$
  (lower curve); the computed slopes are, respectively, $a=3.0$ and
  $a=4.1$.}
\end{figure}

\subsection{$\chi$ for the \av\ model}

\noindent
We studied several \av\ models, among them:
\begin{tabbing}
  \hspace{20mm}\=(a)\qquad $\alpha_i=\pm1$ \hspace{40mm}
               \=(b)\qquad $\alpha_i=1/2\pm1$ \\
               \>(c)\qquad $\alpha_i=1\pm1/2$
               \>(d)\qquad $\alpha_i=1\pm1/3$ \ , \\
\end{tabbing}

\vspace{-4mm}
\noindent
randomly with equal probability. (We also considered models with the
$\alpha_i$ modulated by a few Fourier components; results are similar,
but a little less clear.) Different models do not behave
identically: $\chi$ always follows asymptotically a power law
$\chi\sim\eps^a$, but the exponent $a$ does depend on the model. More
precisely:

\begin{figure}
\centerline{\epsfig{figure=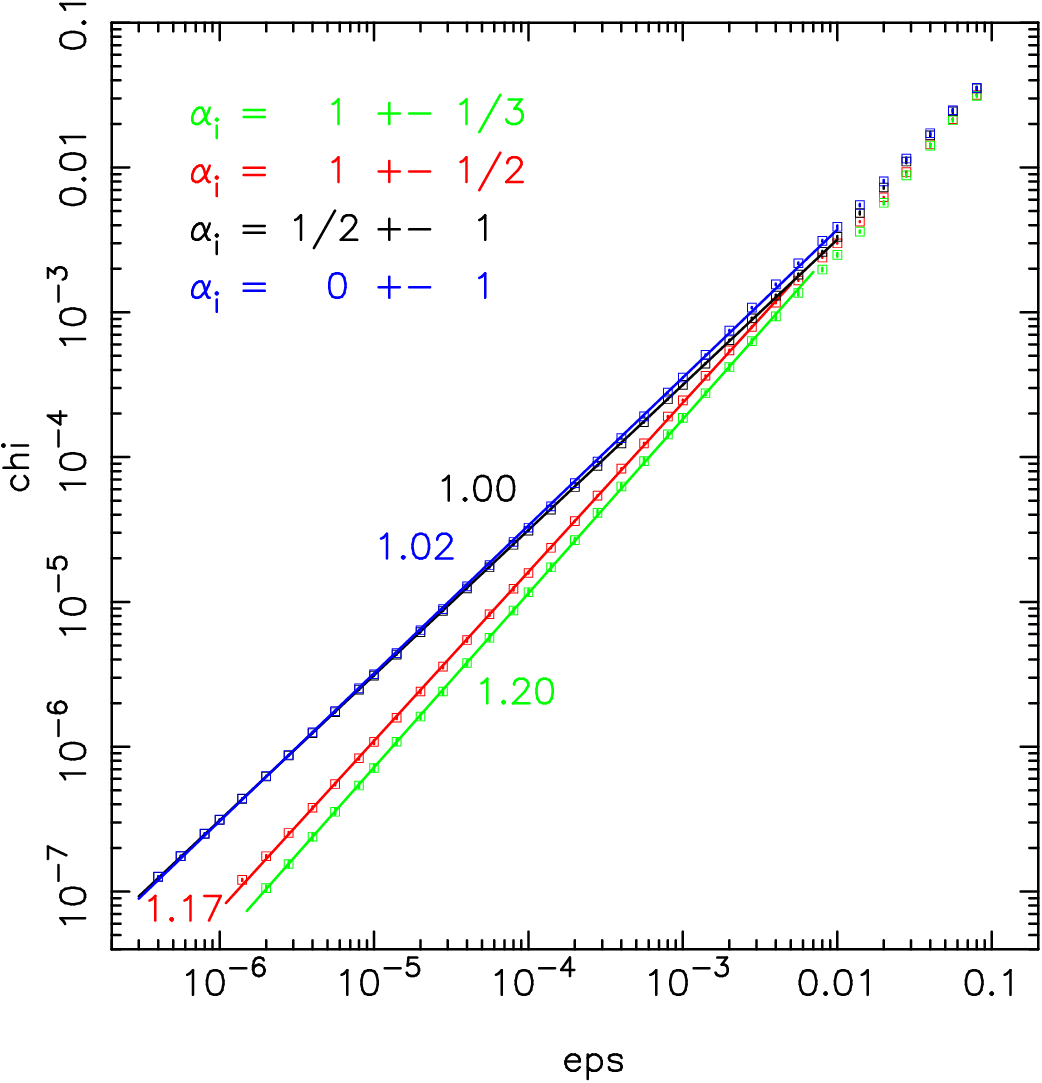,width=80mm}}
\Caption{F-12}{$\chi(N,\eps)$ as function of $\eps$ for the
  \av\ models (a)-(d), $N=8192$. Lower curves: models (a) and (b);
  upper curves (almost superimposed): models (c) and (d).
  Computed slopes, in the order, $a=1.00$, $1.02$, $1.17$, $1.19$. }
\end{figure}

\begin{itemize}
\item[-- ] For models (a) and (b) of the above list, we found
  $a=1$. Such a value fits the rule \calleq{a-lh} of the linear
  hierarchy, actually extending it to $s=1$, and suggests models (a)
  and (b) are part of it, Toda being apparently too far and not
  influencing the low $\eps$ dynamics. See the two upper lines of
  figure \ref{F-12} (hardly distinguishable, for the data almost
  superimpose).

\item[-- ] Models (c) and (d) instead have a larger exponent
$a\simeq1.2$: a somehow misterious value, as all values
of $a$ for the Toda hierarchy are; see the two lower lines of figure
\ref{F-12}. The difference with respect to (a) and (b) suggests the
models feel the presence of Toda, and are part of its
hierarchy, although, in lack of any interpretation of the exponents of the
Toda hierarchy, it is hard to draw any conclusion.

\end{itemize}

\noindent
Our feeling is that the \av\ models feel the proximity of the linear
model, when the variance of $\alpha_i$, in some sense to be
better understood, dominates on the average; by the way, in such a case the
average seems to be irrelevant, as is shown by the exact superposition
of data for models (a) and (b).  If instead the average gets
important, the presence of Toda gets relevant. In lack of any
theoretical basis, we did not further
investigate this point.

\section{Theoretical items}

As remarked in the Introduction, there exists an approximate theory, based
on statistical and geometric considerations, adapted to weakly
nonlinear systems with many degrees of freedom like FPU, which leads
to an analytic estimate of the maximal Lyapunov exponent
\cite{CLP95}\cite{CCP96}\cite{LLPR08}.
The theory was tested on the \pb\ FPU model (as well as in the
so-called $\phi^4$ model and in a chain of weakly coupled rotators)
and the agreement with numerical data turned out to be excellent,
both for small and for large specific energy.

The underlying idea is that, for systems like FPU, the exponential
separation of nearby trajectories is not produced (as in Anosov
systems) by a negative curvature of the manifold where the motion occurs,
i.e.~of the constant energy surface, rather it is a phenomenon of
parametric instability due to the fluctuations of the 
curvature along trajectories, in a convenient natural
metrics. Schematically, the theory proceeds as follows:

\begin{itemize}

\item[i. ] Within certain assumptions (the most important one is that,
in the chosen coordinates, the oscillations or rotations
of the tangent vector $\xi$ do not play an important role),
the problem is reported to a
Hill equation for a single oscillator, namely \beq{hill} \ddot \psi +
\Omg(t)\psi = 0\ , \eeq where however $\Omg(t)$ is not an assigned
function of time but a stochastic process; $\psi(t)$ represents any of
the components of $\xi$.  The process is assumed to be Gaussian and
$\delta$--correlated on a convenient time scale $\tau$. Its average
$\Omg_0$ and variance $\sigma$ (after a nontrivial geometrical
analysis) are found to be
$$
  \Omg_0=\frac1N\media{\Delta \cV}\ ,\qquad
  \sigma^2=\frac1N\big[\media{(\Delta \cV)^2}-\media{\Delta V}^2\big]\ ,
$$
where $\Delta \cV$ is the Laplacian of the potential energy 
$$
\cV(q)=\sum_{i=1}^NV(r_i) \ , \qquad r_i=q_i-q_{i-1} \ ,
$$

and $\media{\,.\,}$ denotes microcanonical averaging. Concerning the
correlation time $\tau$, dimensional arguments lead to two candidates
\cite{CLP95}:
$$
  \tau_1=\sqrt{\frac{2}{\Omg_0}}\ ,\qquad
    \tau_2=\sqrt{\frac{\Omg_0}{2\sigma^2}}\ .
$$
The choice of the authors, leading to good results, is
\beq{tau12}
\tau^{-1}=(\tau_1^{-1}+\tau_2^{-1})\ .
\eeq
The idea underlying the composition rule \calleq{tau12} is that, if
two very different correlation times enter the dynamics, $\chi$ is
sensitive to the shorter one. In the asymptotics $\eps\to 0$ one has
$\sigma\ll\Omg_0$, and thus $\tau=\tau_1$. The numerical factors 
looks a little arbitrary; on the other hand, any dimensional estimate
includes by itself an arbitrary factor.

\item[ii. ] Once the process is defined, a result by Van Kampen \cite{VK76}
applies, leading to
$$
\chi(\Omg_0,\sigma,\tau)=\frac12\Big(\Lambda-\frac{4\Omg_0}{3\Lambda}\Big)\ ,
$$
where
$$
\Lambda=\Big[2\tau\sigma^2
    +\sqrt{(4\Omg_0/3)^3+(2\tau\sigma^2)^2}\,\Big]^{1/3} \ .
$$
\end{itemize}

\noindent
Asymptotically for small $\eps$ one has $\sigma\ll\Omg_0$ and correspondingly
\beq{asintotica}
\Omg_0\to 2\ , \qquad \tau\to 1\ , \qquad
\chi\simeq \frac18\,\sigma^2\ .
\eeq
So, everything is reported to estimating the variance of
$$
\Delta \cV=2\sum_{i=1}^NV''(r_i)-V''(r_1)-V''(r_N)\ ;
$$
the last two terms are due to the choice of fixed ends, and are
negligible for large $N$.  Quite clearly, for the potentials we are
dealing with,
\beq{Delta-V}
\Delta\cV=2N+4\sum_{i=1}^N\alpha_ir_i+6\beta\sum_{i=1}^Nr_i^2
   +8\gamma \sum_{i=1}^N r_i^3+10\delta \sum_{i=1}^N r_i^4 \ ,
\eeq
where the possibility of site depending $\alpha_i$ has been taken
into account.

Let us forget for a moment the constants and look only at the
dependence of $\sigma$, and thus of $\chi$, on $\eps$. The computation gets
straightforward: since at small $\eps$ the harmonic energy dominates,
it is $r_i\sim \sqrt{\eps}$, and consequently
$$
\chi \sim \sigma^2 \sim \eps^{s}\qquad \hbox{for}\qquad
V(r)=\tfrac12r^2+\cO(r^{2+s})\ , \quad\qquad s\ge 1\ .
$$
This is precisely rule \calleq{a-lh} governing the linear
hierarchy. {\it For $s=1$, however, it is necessary that the constants
  $\alpha_i$ do depend on the site}: if $\alpha$ is constant, as in
the \ab\ model, then, because of the fixed ends, the linear term in
\calleq{Delta-V} exactly vanishes (the same would happen with periodic
boundary conditions, actually whenever $\sum_i r_i$ is not allowed to
fluctuate). For the \ab\ model, the theory predicts $\chi\sim \eps^2$,
precisely as for the \pb\ model.

In fact, the theory predicts the exponent $a=2$ not only for the
\ab\ model, but for all models of the Toda hierarchy, {\it including
Toda itself}. This is not surprising: indeed the theory captures the
nonlinearity, rather than the lack of integrability, and cannot be
sensitive to the absolutely peculiar properties of Toda and of the
models close to it.

Computing the exact asymptotics of $\sigma$ in $\eps$, including the
constants, can be done for all models of the linear hierarchy. The
computation we did is sketched in the Appendix; it takes into account
the possibility that not only $\alpha$, but also $\beta$, $\gamma$ and
$\delta$ depend on the site.  Denoting by $\overline{(\cdot)}$ and by
$\sigma^2_{(\cdot)}$, respectively, the arithmetic mean and the
variance of the variables $(\cdot)$, the asymptotic result for $\chi$
turns out to be
\beq{chitab}
\chi(\eps) \simeq \left\{
\begin{array}{ll}
2\sigma_\alpha^2 \,\eps   &\qquad \hbox{\av} \\
\\
\tfrac92\Big(\overline{\beta^2}+\sigma_\beta^2\Big)\,\eps^2
     & \qquad \hbox{\pb} \\
     \\
48\Big(\overline{\gamma^2}+\tfrac32\sigma_\gamma^2\Big)\,\eps^3
     & \qquad \hbox{\gd} \\
     \\
750\Big(\overline{\delta^2}+\tfrac35\sigma_\delta^2\Big)\,\eps^4
     & \qquad \hbox{\pd}
\end{array}
\right.\ .
\eeq
For the \pb\ model, this asymptotic result
reproduces well the asymptotic
line reported in \cite{CLP95}\cite{CCP96}\cite{LLPR08}. Such a line
fits well our numerical data for $N=1024$ or $2048$, while for larger
$N$ the numerical data lie slightly above it; see figure \ref{F-13},
upper-left panel. For the other models of the linear hierarchy, the
situation is as follows:

\begin{figure}
\epsfig{figure=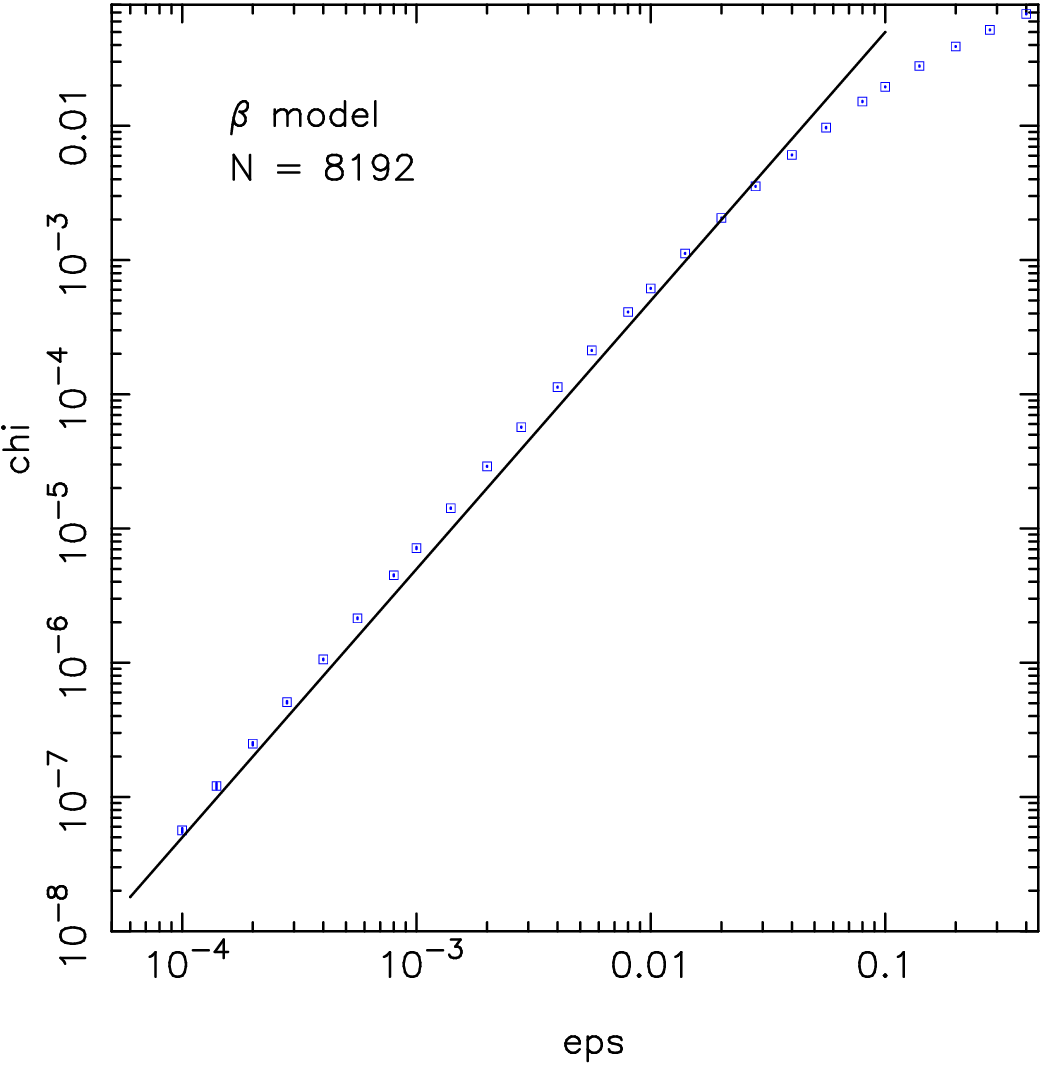,width=80mm}
\hspace{3mm}
\epsfig{figure=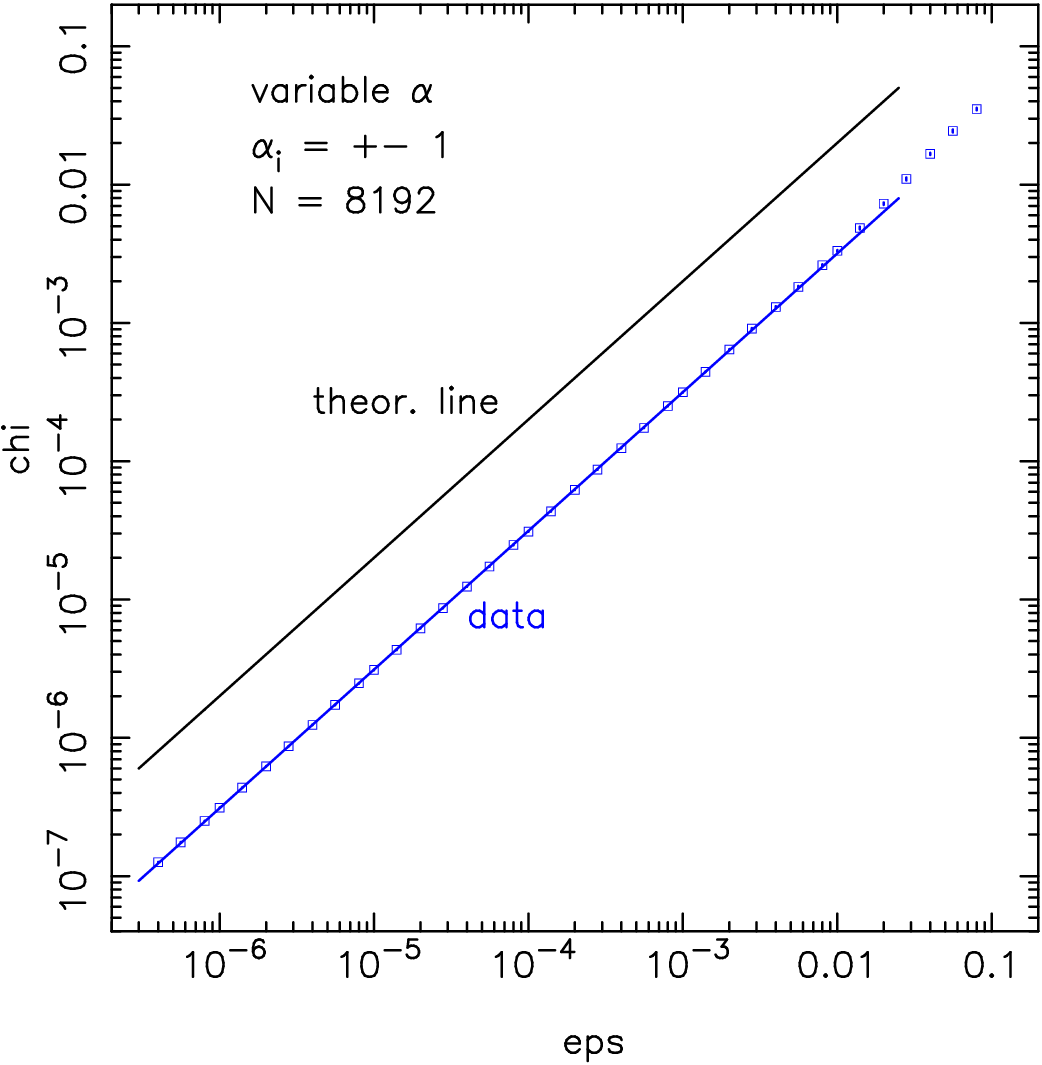,width=80mm}
\vspace{1mm}
\epsfig{figure=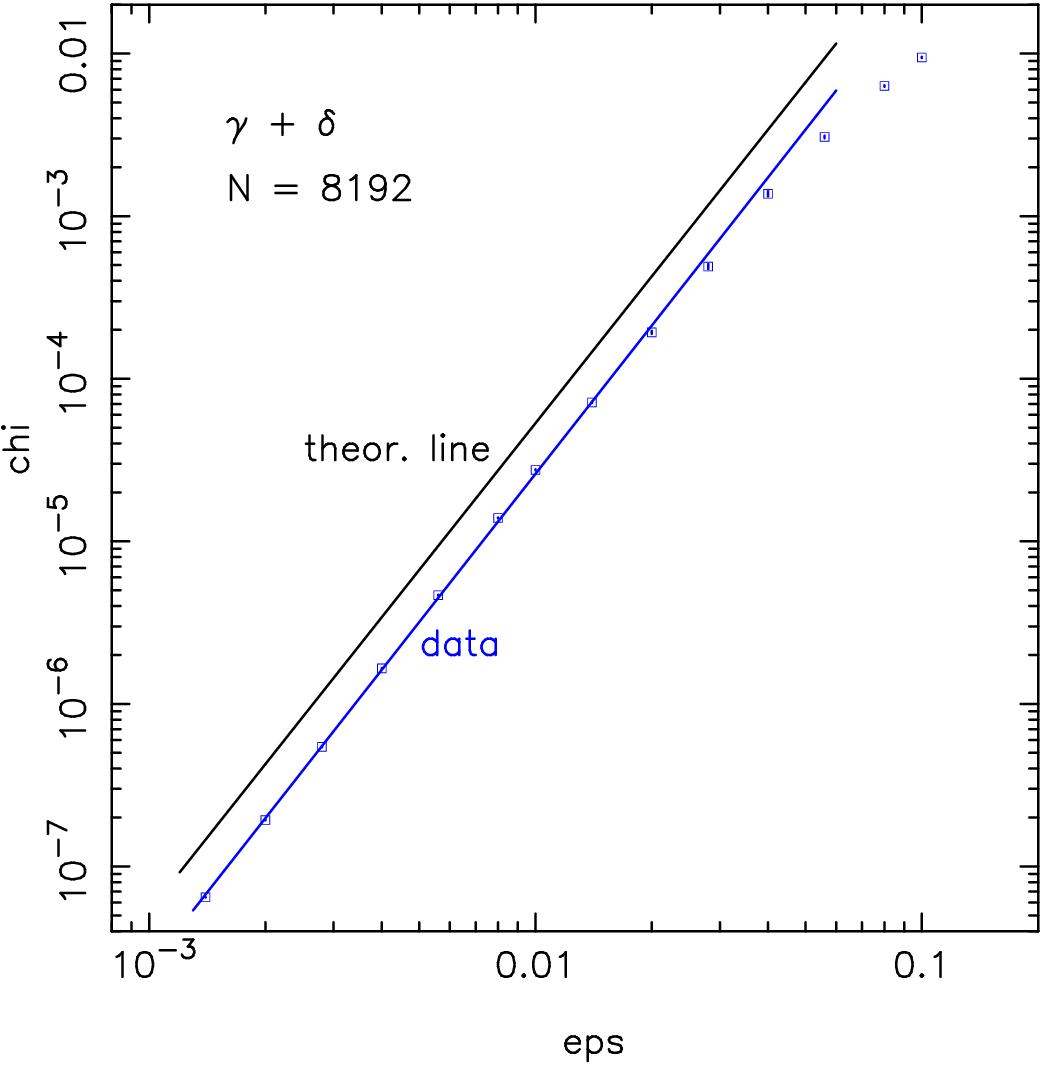,width=80mm}
\hspace{3mm}
\epsfig{figure=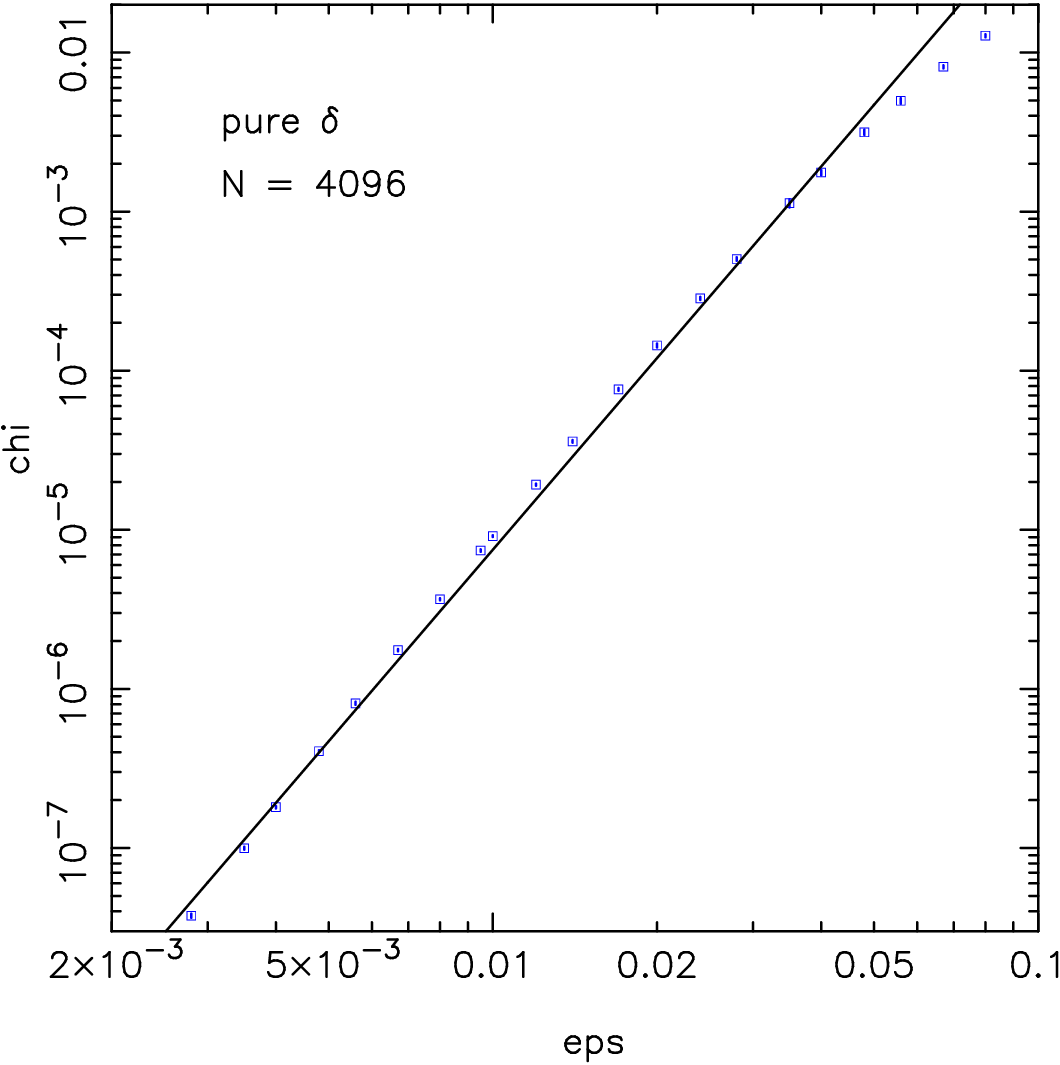,width=80mm}
\Caption{F-13}{A comparison between the theoretical lines and the
  numerical data. Upper left: the \pb\ model; upper right: the
  \av\ model (a); lower left: the \gd\ model; lower right: the
  \pd\ model. $N$ it the largest computed, namely $8192$ for the first
  three models, $4096$ for the last one.}
\end{figure}

\begin{itemize}
\item[-- ] The \av\ models (a) and (b): the theoretical line lies somehow
  above the data, approximately by a factor $6$.

\item[-- ] The \gd\ model: the theoretical line lies
  above the data, approximately by a factor $2$.

\item[-- ] The \pd\ model: the agreement looks quite good.

\end{itemize}

\noindent
See the other panels of figure \ref{F-13}. We are not surprised that a
factor of order one is needed, due to the necessarily rough definition
of $\tau$, and we consider the theory developed in the above quoted
references to be remarkably correct and successful, 
for all models of the linear hierarchy. 

\section{Concluding remarks}

\noindent
{\it A. Numerical errors?} \quad
Numerical integration of differential equations is never safe. The
procedure we used, mentioned at the beginning of Section 2, is however
among the most standard ones, both for trajectories and for the
computation of Lyapunov exponents. 

Concerning trajectories, let us say we successfully used fourth-order
leap-frog in several papers, also investigating theoretically the behavior
of this algorithm in connection with FPU-like models \cite{BP11-2}. Among
the crucial tests, the algorithm is observed to preserve well not
only energy, but also the Toda constants of motion, when used to
integrate Toda (it is obviously essential that the diffusion transversal to
Toda tori is produced by the dynamics, not by the algorithm).

Concerning specifically the computation of Lyapunov exponents, the
main risk, as far as we understand, is overestimation: consider the
case in which a small or
vanishing exponent is due to the fact that, in the metrics at hand, a
relatively large expansion (easy example: local hyperbolicity,
tangent vector pointing in an expanding direction) is followed, after
some time $T$, by a compensating contraction (local hyperbolicity, evolved
tangent vector pointing now in a contracting
direction). Then, if $T$ is large, there is a chance such a
correlation is lost, compensation is partially missed and expansion prevails.

A good test is then computing $\chi$ for Toda, for which compensation
is exact, and looking how close to
zero is the computed value; intentionally, the integration step $\ts$
can be taken large, so as to emphasize the effect on $\chi$ of a
rough algorithm. Results are in figure \ref{F-14}. The left panel shows
$\chib(t)$ vs. $t$ for the Toda model with $N=1024$, at
$\eps=8\times10^{-4}$, $\ts$ ranging from $0.05$ to $0.40$; for each
$\ts$, the {\it plateau} reached by $\chib(t)$ provides an estimate
of the possible error in $\chi$ in models close to Toda. The figure
clearly resembles figure \ref{F-3}, as is not surpring since the
essential effect of the integration algorithm is that of introducing a
perturbation \cite{BG94} \cite{H94} \cite{HLW02},
and correspondingly, the integrated system can be
considered a further member of the Toda hierarchy. The computed value
of $\chi$ (the error) also depends on $\eps$, as is shown in the right
panel of the figure, where $\chi$ is reported vs. $\eps$ for fixed
$\ts=0.24$ (a quite large value, so as to have large enough computable
errors). As is remarkable, here too we find a well defined power
law $\eps^a$, with $a\simeq1.60$. The value of $a$, almost identical
to $a$ of the \ab\ model at the same $N$, suggests that the effective
perturbation on Toda introduced by the algorithm is quartic. We did not
further investigate this point, although it could be done along the
lines of \cite{BP11-2}.

Practically, in our numerical computations we always took care that
the error due to the algorithm was much smaller than the measured value of
$\chi$. In most cases, $\ts=0.1$ (the typical value we used with
confidence, in the past, with fourth order leap-frog) was fairly
enough; only for the \gt\ and the \pd\ models, at small $\eps$,
it was necessary to use $\ts=0.05$.

\begin{figure}
\epsfig{figure=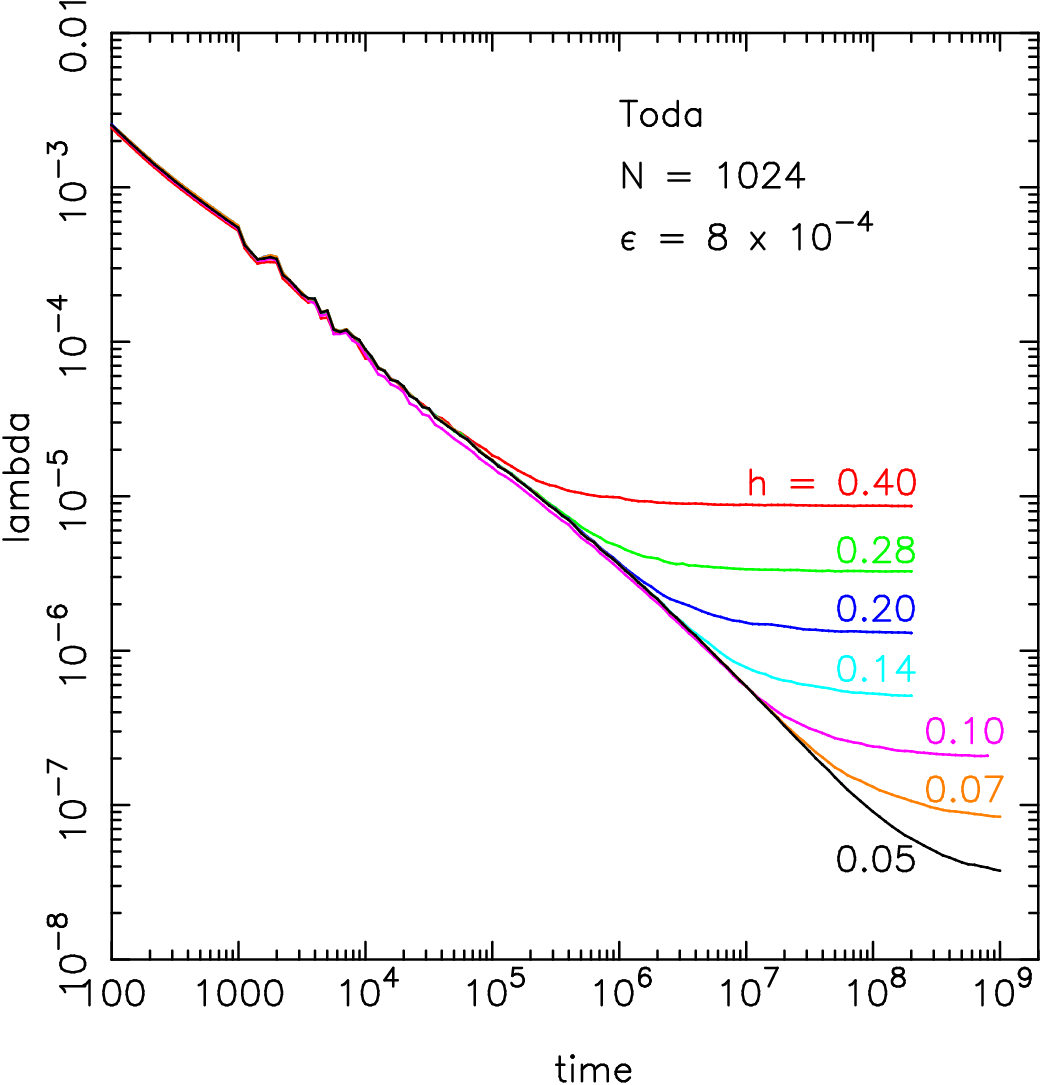,width=80mm}
\hspace{3mm}
\epsfig{figure=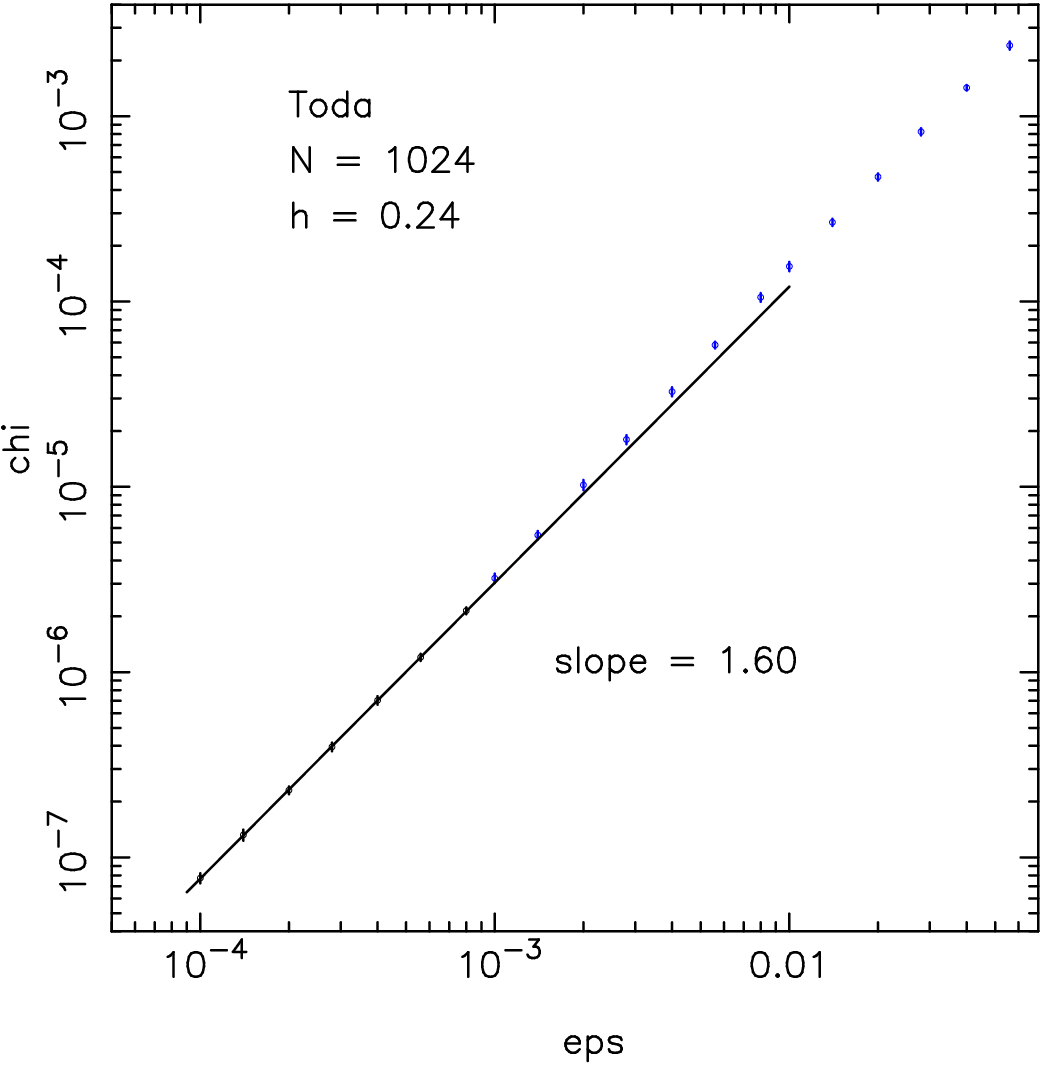,width=80mm}
\Caption{F-14}{Left: $\chib(t)$ vs. $t$ for the Toda model with
  $N=1024$, at $\eps=8\times10^{-4}$, for several choices of the
  integration step $\ts$ between $0.05$ and $0.4$. Right: the
  dependence of the limit value $\chi$ on $\eps$, for fixed (large)
  $\ts=0.24$. }
\end{figure}

\bl

\noindent
{\it B. Anomalous exponents for the Toda hierarchy?}\quad
The power laws we produced look rather precise. For the linear
hierarchy the exponents are ``simple'' and theoretically clear,
while for the Toda hierarchy
they are not. For the \ab\ model we observed a slow asymptotics in $N$
leading to $a=3/2$,
while for the other models of the hierarchy simple
values were not observed. As far as we know, all power laws met in
FPU, numerically or theoretically, in the weak stochasticity regime,
have exponents either integer, or half-integer, or with denominator
$4$. The exponents we found for the \bt\ and the \gt\ models, as well
as for the \av\ models (c) and (d), look then anomalous and are hard
to understand. Some theory, similar and as powerful as the one
reported in \cite{CLP95}\cite{CCP96}\cite{LLPR08}, but adapted to the
neighborhood of the Toda model, seems to be necessary,
but developing it looks a rather difficult task.

\vspace{10mm}
   
 \section*{The Appendix: analytical details}
\renewcommand{\thesection}{A}
\numberwithin{equation}{section}
\setcounter{equation}{0}

In this Appendix we report some details concerning the analytical
estimate of the Lyapunov exponents for the linear hierarchy.  As
stressed above, the theory is essentially that developed in
\cite{CLP95,CCP96} for the pure $\beta$-model: we just extended it to
the other models of the hierarchy, limiting ourselves to
compute the asymptotic behavior of $\chi$ as $N\to\infty$ and
$\eps\ll1$. The symbol '$\simeq$', in the sequel, means asymptotic
equality in this limit. As pointed out in Section 3, according to the
theory the maximal Lyapunov exponent is given by
\beq{chiform}
\chi(\eps)\simeq\frac{\sigma^2}{8}
=\frac{
\media{(\Delta \cV)^2}-\media{\Delta \cV}^2}{8N}\ ,
\eeq
where $\media{\,\cdot\,}$ denotes the microcanonical average 
on the constant energy surface $H(p,q)=N\eps$. 
The right hand side of \calleq{chiform} has to be computed to leading order in
$\eps$. 
We admit the possibility that the constants, at any order of nonlinearity,
depend on the site, that is
\beq{cV}
\cV=\sum_{i=1}^NV_i(r_i)=\sum_{i=1}^N\left(\frac{r_i^2}{2}
+\alpha_i\frac{r_i^3}{3}+\beta_i\frac{r_i^4}{4}
+\gamma_i\frac{r_i^5}{5}+\delta_i\frac{r_i^6}{6}\right)\ ,
\eeq
only assuming that at any site, the pair potential $V_i(r)$ displays
a single minimum at $r=0$ and diverges to $+\infty$ for
$r\to\pm\infty$ ($\delta_i>0$, or
$\delta_i=\gamma_i=0$ and $\beta_i>0$).
The corresponding expression of the Laplacian is
\beq{DelcV}
\Delta\cV=2N+4\sum_{i=1}^N\left(\alpha_ir_i+6\sum_{i=1}^N \beta_ir_i^2
   +8\sum_{i=1}^N\gamma_ir_i^3+10\sum_{i=1}^N\delta_ir_i^4
   \right)\ .
\eeq
As pointed out in
\cite{CLP95,CCP96}, the microcanonical quadratic
fluctuation of $\Delta\cV$, appearing in the numerator on the right
hand side of \calleq{chiform}, can be computed in terms of suitable
\emph{canonical} averages, according to the theory of Lebowitz, Percus
and Verlet (LPV)~\cite{LPV67}. For $\eps\ll1$
the LPV formula reads, to leading order in $\eps$,
\beq{LPVDelcV}
\frac{\media{(\Delta\cV)^2}-\media{\Delta\cV}^2}{N} \simeq
\frac{\media{(\Delta\cV)^2}_c-\media{\Delta\cV}_c^2}{N}-
\eps^2\left(\frac{\media{\Delta\cV}_c'}{N}\right)^2\ ,
\eeq
where $\media{\,\cdot\,}_c$ denotes the canonical average at
temperature $\eps$, $\media{\Delta\cV}_c'=d\media{\Delta\cV}_c/d\eps$,
and it is understood that all quantities on the right hand side are
computed to leading order in $\eps$. Using the simple scaling argument
$r_i\sim\sqrt{\eps}$, valid when $\eps\ll1$, from expression
\calleq{DelcV} one easily realizes that, according to the LPV formula
\calleq{LPVDelcV}, the microcanonical and the canonical quadratic
fluctuation of $\Delta\cV$ differ, to leading order, only for leading
nonlinearity of even order, that is for the $\beta$-models ($\alpha_i=0$),
and the $\delta$-models ($\alpha_i=\beta_i=\gamma_i=0$).

Up to this point we did not introduce any novelty, apart from extending
formulas to all models of the linear hierarchy. We now proceed by
computing the statistical averages for the different models; for the
sake of consistency, we shall take into account the constraint
$\sum_{i=1}^Nr_i=0$, which holds both for fixed
ends and for periodic boundary conditions, whereas in~\cite{CLP95,CCP96}
the analytic computations are made in the simpler case of free
ends. Such a constraint, as already mentioned in Section 3, plays a
subtle role.

All canonical averages on the right hand side of
\calleq{LPVDelcV} are clearly sums, with easily computed coefficients,
of moments of the form
$$
\media{r_j^mr_k^n}_c=\frac{\int r_j^{m}r_k^{n}\ 
e^{-\frac{1}{\eps}\cV(r_1,\dots,r_N)}
\delta(\sum_ir_i)\ dr_1\dots dr_N}{
\int e^{-\frac{1}{\eps}\cV(r_1,\dots,r_N)}
\delta(\sum_ir_i)\ dr_1\dots dr_N}\ ,
$$
where $j,k$ are any indices from $1$ to $N$, while $m,n$ are exponents
between $0$ and $4$ and $\delta(x)$ denotes the Dirac function (expressing the
constraint $\sum_{i=1}^Nr_i=0$), all integrals running from $-\infty$ to
$\infty$.  By representing the Dirac function in the integral form
$\delta(x)=\frac{1}{2\pi}\int e^{iyx}dy$, using the expression
\calleq{cV} of the potential energy, one gets \beq{momjk2}
\media{r_j^mr_k^n}_c=(-\imath)^{m+n}\ \frac{\int
  \Big[\prod_{\substack{i=1 \\ i\neq j,k}}^Nf_i(y)\Big]
  f_j^{(m)}(y)f_k^{(n)}(y)\ dy}{ \int \prod_{i=1}^Nf_i(y)\ dy}\ , \eeq
where $\imath$ denotes the imaginary unit, $f_i$ is defined as
\beq{f-i}
f_i(y)=\int e^{-\frac{1}{\eps}V_i(r)+\imath yr}\ dr \ ,
\eeq
and $f_i^{(s)}$ denotes the $s$-derivative of $f_i$. 
Now, under the hypothesis made above that
$V_i(r)$ has an unique minimum at $r=0$, one easily shows, by the Laplace
method (see for example \cite{AF03}), that for small $\eps$ it is 
$$
f_i(y)\simeq \sqrt{2\pi\eps}\ e^{-\eps\frac{y^2}{2}} \ ,
$$
independent of $i$ (this is not surprising since in \calleq{f-i},
for small $\eps$, only the $i$-independent harmonic part of $V_i$ 
contributes to the integral). In such a way, the average
\calleq{momjk2} assumes the simple asymptotic form
\beq{momjk3}
\media{r_j^mr_k^n}_c=(-\imath)^{m+n}\
\frac{\int e^{-N\eps\frac{y^2}{2}}\,P_{m+n}(y)\ dy}
{\int e^{-N\eps\frac{y^2}{2}}\ dy}\ ,
\eeq
where $P_{m+n}(y)$ is a suitable polynomial of degree $m+n$; moreover,
it turns out that if $m+n$ is even (odd), then $P_{m+n}(y)$ contains
only even (odd) terms, so in particular, for odd $m+n$, the right hand
side of \calleq{momjk3} vanishes.

In this way one reduces all the analytic work to computing trivial
Gaussian moments.  Taking a bit carefully into account the dependence
on $N$ and summing up all terms, yields the asymptotic estimates
\calleq{chitab}.

\end{document}